\documentclass[11pt]{article}
\usepackage[margin=2.5cm]{geometry}
\usepackage{pifont}
\usepackage{tikz-cd}

\usepackage{caption}
\usepackage{microtype}
\usepackage{xcolor}
\definecolor{Red}{rgb}{1,0.0.25,0.25}
\definecolor{Green}{rgb}{0.25,0.75,0.25}
\definecolor{Blue}{rgb}{0.1,0.5,1}
\usepackage[normalem]{ulem}
\usepackage{amsmath}
\usepackage{mathtools}
\usepackage{amssymb}
\usepackage{amsthm}
\usepackage{thmtools}
\usepackage{graphicx}
\graphicspath{ {images/} }
\usepackage{commath}
\usepackage{setspace}
\usepackage{float}
\usepackage[final]{pdfpages}
\usepackage{bm}
\PassOptionsToPackage{hyphens}{url}
\usepackage{xparse}
\usepackage[english]{babel}
\usepackage[colorinlistoftodos]{todonotes}
\usepackage{csvsimple}
\usepackage{mwe}
\usepackage{float}
\usepackage[version=3]{mhchem}
\usepackage{bm}
\usepackage{mathrsfs}
\usepackage{tcolorbox}
\usepackage[thinc]{esdiff}
\usepackage{lscape}
\usepackage{tabu}
\usepackage{multirow}
\usepackage{lipsum}
\usepackage{subcaption}
\usepackage{graphicx}
\usepackage{hyperref}
\hypersetup{
  colorlinks = true,
  linkcolor =blue,
  anchorcolor = red,
  citecolor = blue,
  urlcolor = blue
}
\usepackage{cleveref}
\usepackage{MnSymbol}
\usepackage{tablefootnote}
\usepackage{multirow}
\usepackage{array}
\numberwithin{equation}{section}

\newtheorem{definition}{Definition}
 
\newtheorem{remark}{Remark}

\newtheorem{lemma}{Lemma}

\newtheorem{theorem}{Theorem}
\newtheorem{assumption}{Assumption}

 \newcommand{\rank}[1]{\operatorname{rank}(#1)}
 
\begin{document}
 \title{Coordinate-independent model reductions of
chemical reaction networks based on geometric singular perturbation theory}
\author{Timothy Earl Figueroa Lapuz\thanks{School of Mathematics and Statistics, The University of Sydney, Camperdown NSW 2006, Australia} \and Martin Wechselberger\footnotemark[1]}
  
\maketitle
\begin{abstract}
The quasi-steady-state approximation (QSSA) is a standard technique for reducing the complexity of chemical reaction networks (CRNs). The validity of any QSSA-based model is restricted to specific parameter regimes. {Selecting the appropriate reduction is not always straightforward. At times, QSSAs are misused outside of their validity regions and, even when a particular QSSA is considered valid in a given parameter regime, other QSSAs may be simultaneously valid, creating ambiguity}.

Here, we employ a more powerful alternative: a constructive model reduction framework based on coordinate-independent geometric singular perturbation theory (ci-GSPT) and the parametrization method. A key advantage of this {approach} is its ability to derive reduced models independent of a clear timescale separation in the variables {for a specific parameter configuration}. 

We demonstrate our approach on two benchmark systems. For the Michaelis-Menten (MM) reaction, we show that the framework provides a {systematic approach by exploring parameter configurations across three orders of magnitude: asymptotically large, small, and `order one'. A consequence of this systematic analysis is a geometric classification that categorizes the resulting model reductions and provides a point of comparison between our approach and common QSSA variants in the literature.
For the more complex Kim-Forger model, we show that {this approach} successfully produces a reduction without the need for a coordinate transformation, showcasing its applicability to larger systems.}

\textbf{Relevance to Life Sciences.}
CRNs are essential for modeling biological processes, from enzyme kinetics (e.g., the MM reaction scheme) to rhythmic regulation of cell functions (e.g., the Kim-Forger model of the circadian clock). Reducing the complexity of these networks is crucial for making mathematical analysis and numerical simulation computationally tractable. While methods like the QSSA provide such reductions, they require the selection of a suitable model reduction based on the system's parameters. Our {chosen approach} circumvents this problem by calculating appropriate model reductions for a specific parameter configuration. This ensures that the reduction is mathematically sound and removes the guesswork from the modeling process. {A catalog of valid MM model reductions arising from our analysis is provided for a variety of parameter configurations}.

\textbf{Mathematical Content.}
Our approach is built on the rigorous foundation of Tikhonov-Fenichel theory for singularly perturbed systems. We utilize two geometric tools that do not require an explicit coordinate transformation into slow and fast variables: (i) ci-GSPT, which allows us to geometrically identify {the critical/slow manifold and compute the reduced dynamics of the system,} and (ii) the parametrization method, which is used to systematically compute higher-order, more accurate approximations of both the invariant slow manifold and the slow vector field defined on it. Together, these tools provide a powerful and rigorous {framework} for deriving model reductions by geometrically identifying a system's invariant slow manifold and its corresponding slow flow, extending the analysis {beyond methods based purely on timescale separation in the variables.}
\end{abstract}

\section{Introduction}
\label{sec:intro}
The {\em reversible Michaelis-Menten (MM) reaction scheme}
\begin{align}
    \ce{S + E <=>[k_1][k_{-1}] C <=>[k_2][k_{-2}] P + E}. \label{crn2}
\end{align}
is one of the most well-known reactions in biochemistry. Here, $S$ refers to a `Substrate', $E$ to an `Enzyme', $P$ to a `Product', $C$ to an intermediate `Complex', and $k_1,k_{-1},k_2 >0$ and $k_{-2} \geq 0$ are reaction rate constants. The {\em irreversible MM reaction scheme}, obtained by setting $k_{-2} = 0$ in \eqref{crn2}, i.e.,
\begin{align}
    \ce{S + E <=>[k_1][k_{-1}] C ->[k_2] P + E} \label{crn}
\end{align}
is even more familiar. Pioneering works by Henri \cite{henri} and Michaelis and Menten \cite{michaelismenten} in the early 1900s derived approximations for the substrate depletion rate. A century later, obtaining long-term approximations for the species concentrations in {\em chemical reaction networks (CRNs)} remains a central challenge. These {\em model reductions} are highly sought after because they reduce the dimension of the governing ordinary differential equations (ODEs), making mathematical analysis and numerical simulation more tractable.

A typical analysis of a CRN begins by applying the {\em law of mass-action (LMA)} to derive a system of ODEs, followed by an exact reduction using conserved quantities. One then proceeds to find an approximate model reduction via, for example, a {\em quasi-steady-state approximation (QSSA)}. For the MM reaction schemes \eqref{crn2} and \eqref{crn}, two common QSSAs are invoked: the standard QSSA (sQSSA) and the total QSSA (tQSSA) \cite{briggshaldane,hta,segelslemrod,borghans,tzafriri}. These approximations have corresponding validity conditions. {A significant gap in understanding emerges due to QSSAs being misused outside of their validity regime \cite{kimtyson}. Further ambiguity emerges due to overlapping validity regimes where more than one QSSA is valid \cite{accuracy}. Two crucial questions arise:} (a) which QSSA should be used in a given parameter region, and (b) can this approach be reliably applied to larger CRNs?

In this paper, we approach these questions differently, using a modern quantitative framework based on {\em coordinate-independent geometric singular perturbation theory} (ci-GSPT), complemented by the {{\em parametrization method} \cite{goekewalcher,feliukruffwalcher,wechselberger2020,coulletspiegel,cabre20031,cabre20032,cabre2005,multiple}}. A crucial feature of this  {integrated approach  \cite{multiple}} is that it is coordinate-independent; unlike classical methods, no coordinate transformations are required to put the system into a standard singularly perturbed form\footnote{A singularly perturbed system is in {\em standard form}, if there is a clear timescale separation for the dynamics of the (slow and fast) variables. In general, transforming ODEs to standard form is non-trivial and is only guaranteed locally; see e.g. \cite{fenichel,wechselberger2020} for details.}, a feature that is often desirable in applications. This allows us to calculate model reductions for a specific parameter configuration, circumventing the need {to select appropriate QSSAs and thereby avoiding misuse outside its validity regime and choosing between competing QSSAs}. The parametrization method is an iterative procedure for calculating higher-order terms of the invariant slow manifold and its slow flow. This is especially useful when the leading-order approximation of the slow vector field vanishes -- which we observe in many of our model reductions -- or when the important system dynamics involve more than two distinct timescales.

{In particular, we systematically explore parameter configurations of the MM reaction network across three orders of magnitude: asymptotically large, small, and ‘order one’, illustrating the applicability of our approach and demonstrating that selecting a valid
QSSA \textit{a priori} is no longer required. As a consequence of this analysis, we present a classification that categorizes a given parameter configuration based on the underlying geometry in the singular limit. This is then compared to the QSSA taxonomy in the literature. Finally, the applicability of our chosen approach to larger systems is demonstrated via calculating a model reduction for the Kim-Forger (KF) reaction scheme, which is a circadian clock model as detailed in, e.g., \cite{kimforger,kim2021}.}

We do note that there are powerful reduction techniques out there, including the {\em computational singular perturbation (CSP)} method\footnote{This is an algorithmic method for the identification of the slow manifold and its slow flow; it shares similarities with the parametrization method.} {\cite{lamgoussis1989,lamgoussis1994,zagariskaperkaper,computational}}. In particular, the CSP method has been employed to calculate a reduction for the MM reaction scheme that covers `almost all' of parameter space and is similarly applicable to larger CRNs \cite{mmvalid,algorithmic2023}. Our approach remains distinct: instead of deriving a single, overarching reduction with a corresponding validity condition, we compute the {specific, appropriate reductions for a given parameter configuration.}

The remainder of this paper is structured as follows. Section \ref{sec:reaction} discusses the pre-analysis steps for CRNs. Section \ref{sec:GSPT} provides a summary of {the tools we use}. Section \ref{MM_class} discusses a geometric classification for the MM reaction scheme, and Section \ref{lit_comp_irr} compares this to QSSA validity conditions. We then list the model reductions for the irreversible and reversible MM schemes in Sections \ref{lit_comp_irr} and \ref{sec:rev}, respectively. We finally conclude in Section \ref{sec:conclusion}. 

The main advances of this work are twofold. Mathematically, we {demonstrate a rigorous and systematic approach for reducing CRNs}. Biologically, we provide a definitive catalog of 14 and 25 distinct, relevant model reductions for the irreversible and reversible MM schemes, respectively, and demonstrate our approach on the more complex KF oscillator. 

\section{Pre-processing CRN Models}
\label{sec:reaction}

In this section, we succinctly describe the pre-processing steps applied to a CRN model prior to a ci-GSPT based model reduction. We note that the steps shown here are not new; see e.g. \cite{leeothmer,goekewalcher} for further details and the textbook \cite{foundations} for a thorough exposition of CRNs from a mathematical point-of-view. 

\subsubsection*{The law of mass-action}
Based on the LMA, a CRN
is described by a system of ODEs of the form
\begin{equation}\label{LMA-orig}
    \frac{dX}{dT}=\mathcal S\, V(X)
\end{equation}
where $X\in\mathbb{R}^n_{\geq 0}$ represents the set of chemical species concentrations involved in the reaction network, $\mathcal S$ is the corresponding $n\times m$ stoichiometry matrix with integer entries which encodes the reaction network structure (i.e., the topology of the network), and $V(X)\in\mathbb{R}^m$ is the corresponding reaction vector which encodes the reactions including their rate constants (i.e., the kinetics of the reactions based on LMA).

\subsubsection*{Dimensional analysis}
Non-dimensionalization of system~\eqref{LMA-orig} is accomplished by rescaling the variables $X_i=k_{x_i}x_i$ with reference scales {$k_{x_i} \in\mathbb{R}_{> 0}$}, $i=1,\ldots,n$,
and $T=k_t t$ with reference timescale {$k_t\in\mathbb{R}_{> 0}$} to obtain the dimensionless system
\begin{equation} \label{LMA-dimless}
    \frac{dx}{dt}=\tilde {\mathcal S} \tilde V(x)
\end{equation}
where $x=(x_1,\ldots,x_n)^\top\in\mathbb{R}_{\geq 0}^n$ denotes the dimensionless chemical species vector, $\tilde{\mathcal S}$ the dimensionless $n\times m$ `stoichiometry matrix' with real entries, i.e., not necessarily integers anymore, and $\tilde V(x)\in\mathbb{R}_{\geq 0}^m$ is the dimensionless reaction vector.

An important outcome of such a dimensional analysis is that there is a reduction of the effective free parameters in the corresponding dimensionless system \eqref{LMA-dimless}, which is a consequence of the {\em Buckingham-Pi Theorem} \cite{buckingham}. Furthermore, the resulting dimensionless parameters compare reaction rates and species concentrations relative to each other.

\subsubsection*{Stoichiometric compatibility class (conservation laws)}
\label{objective}

An important observation for CRNs such as \eqref{crn2} is that the network topology encoded in the stoichiometric matrix may include redundant information. If so, then the dynamics of such a CRN is not exploring the whole phase space $\mathbb{R}^n_{\geq 0}$, but is restricted to an affine invariant subspace or, more precisely, a lower-dimensional simplex of the phase space. This leads to a first and {\em exact} model reduction of a CRN to its so-called {\em stoichiometric compatibility class}. 

\begin{assumption}
The $n\times m$ stoichiometry matrix has rank $1\le r < min(n,m)$.
\end{assumption}

For $1\le r < min(n,m)$, there exists a corresponding  $n\times(n-r)$ matrix $L$ such that {$\rank{L} = n-r$} and $L^\top \mathcal{\tilde{S}} = \mathbb{O}_{n-r,m}$. The column vectors of $L$ form a basis of the left null-space for $\mathcal{\tilde{S}}$, and this basis forms the orthogonal complement of the desired affine invariant subspace of the CRN, i.e., $L^\top \frac{dx}{dt} = \mathbb{O}_{n-r}$. The location of this affine subspace 
is determined by the initial condition (IC) $x_0$ of system \eqref{LMA-dimless}.

\begin{definition}
    The $r$-dimensional invariant affine subspace 
    $${\mathcal A}_{\tilde {\mathcal S}}:=\{x\in \mathbb{R}^n_{\geq 0} \, : \, L^\top (x-x_0) = 
    \mathbb{O}_{n-r}\}$$
    is the dimensionless{ \em stoichiometric compatibility class} of system \eqref{LMA-dimless}.
\end{definition}

The dimensionless stoichiometric compatibility class ${\mathcal A}_{\tilde {\mathcal S}}$ defines corresponding $n-r$ dimensionless {\em conservation laws}, i.e., $n-r$ conserved quantities $L^\top x=L^\top x_0$. These $n-r$ (linear) equations can be solved, e.g., as a global graph over $r$ of the given $n$ coordinates. This gives the desired exact $r$-dimensional model reduction to
\begin{equation}\label{LMA-stoich-red}
    \dfrac{dy}{dt} = \mathcal{V}^\top \hat{V}(y) =: \mathcal{F}(y;\mu)
\end{equation}
where $\mathcal{V}^\top$ is a full row-rank $r\times n$ matrix\footnote{There exists a (non-unique) rank factorisation of $\tilde{S} = \mathcal{U}\mathcal{V}^T$ with corresponding full column-rank $n\times r$ matrix $\mathcal{U}$. A left-inverse $\mathcal{U}^L$ is applied in order to reduce system \eqref{LMA-dimless} to system \eqref{LMA-stoich-red}.} and $\hat{V}(y)$ is the transformed dimensionless reaction vector. 
Here, the vector field $\mathcal{F}(y;\mu)\in\mathbb{R}^r$ { is polynomial and} emphasizes its dependence on the dimensionless parameters $\mu\in\mathbb{R}_{\geq 0}^p$, $p\in\mathbb{N}$.
   
\begin{remark}
    The stoichiometric reduction step also applies to the original dimensional model \eqref{LMA-orig} since {$\rank{\mathcal{S}} = \rank{\mathcal{\tilde{S}}} = r < min(n,m)$}.
\end{remark}

\subsection{Pre-processing the reversible MM reaction scheme}

\subsection*{The law of mass-action}
For the CRN \eqref{crn2}, the LMA gives the following system of ODEs,
{
\begin{equation}\label{eq:orig-EK}
    \begin{pmatrix}
\frac{dS}{dT}\\
\frac{dE}{dT}\\
\frac{dC}{dT}\\
\frac{dP}{dT}
\end{pmatrix}
=
\begin{pmatrix}
-1 & +1 & 0 & 0\\
-1 & +1 & +1 & -1\\
+1 & -1 & -1 & +1\\
0 & 0 & +1 & -1
\end{pmatrix}
\,
\begin{pmatrix}
k_1 SE\\
k_{-1} C\\
k_2 C\\
k_{-2}PE
\end{pmatrix}  = \mathcal{S} V(X)
\end{equation}}where $S, C, E$ and $P$ denote the concentrations of chemical species, and the (typical) initial conditions (ICs) are $S(0) = s_0>0$, $E(0) = e_0>0$, $C(0) = 0$ and $P(0) = 0$. The units of the variables and parameters are listed in Table~\ref{units_table}.
\begin{table}[ht]
\centering
\begin{tabular}{ |c|c|c|c| } 
 \hline
Variable/parameter & Unit & Parameter & Unit \\ \hline
S & M (mol$\cdot$L$^{-1}$) & $k_1$ & M$^{-1}$s$^{-1}$ \\ 
E & M & $k_{-1}$ & s$^{-1}$ \\
C & M & $k_2$ & s$^{-1}$ \\
P & M & $k_{-2}$ & M$^{-1}$s$^{-1}$ \\
T & s & & \\
\hline
\end{tabular}
\caption{Units of the variables and parameters for the MM reaction scheme.} \label{units_table}
\end{table}
\subsubsection*{Dimensional analysis}
For system \eqref{eq:orig-EK}, we define the rescaling of the species concentration as follows
$$
S=k_S s,\quad E=k_E e,\quad C=k_C c,\quad P= k_P p,\quad T=k_T t
$$
where $k_S=k_P=s_0$ and $k_C=k_E=e_0$ are reference scales of the chemical species concentrations and $x=(s,e,c,p)^\top$ is the corresponding dimensionless species concentration vector, $k_t= (k_1s_0)^{-1}$ is a reference timescale ({i.e., first forward reaction in the reaction network \eqref{crn2}})
and $t$ the corresponding dimensionless time. 

This choice of reference scales transforms \eqref{eq:orig-EK} into its dimensionless form
\begin{equation}\label{eq:dimless-EK}
\frac{dx}{dt}=
    \begin{pmatrix}
\frac{ds}{dt}\\

\frac{de}{dt}\\

\frac{dc}{dt}\\

\frac{dp}{dt}
\end{pmatrix}
=
\begin{pmatrix}
-\beta & +\beta & 0 & 0\\
-1 & +1 & +1 & -1\\
+1 & -1 & -1 & +1\\
0 & 0 & +\beta & -\beta
\end{pmatrix}
\,
\begin{pmatrix}
se\\
\alpha c\\
\gamma c\\
\delta pe
\end{pmatrix}
= \tilde {\mathcal S} \tilde V(x) 
\end{equation}
with dimensionless initial conditions 
$s(0) = 1$, $e(0) = 1$, $c(0) = 0$ and $p(0) = 0$ and four
non-negative dimensionless parameters defined as
\begin{equation}
\alpha := \dfrac{k_{-1}}{k_1 s_0}, \quad \beta := \dfrac{e_0}{s_0}, \quad \gamma := \dfrac{k_{2}}{k_1 s_0},\quad \delta := \dfrac{k_{-2}}{k_1}\,,
\end{equation}
where
\begin{itemize}
    \item 
$\alpha$ measures the ratio between the first backward and forward reaction:\\ 
\ce{S + E  <=>[k_1][k_{-1}] C},
\item $\beta$ measures the ratio of initial enzyme and substrate concentrations,
   \item 
$\gamma$ measures the ratio between the two consecutive forward reactions:\\
\ce{S + E  ->[k_1] C  ->[k_2] P + E},
\item $\delta$ measures the ratio between the first forward and the second backward reaction: 
\ce{S + E  ->[k_1] C  <-[k_{-2}] P + E}. 
\end{itemize}

\begin{remark}
    Since we have three ratios ($\alpha,\gamma,\delta$) of the four involved reactions in \eqref{crn2}, we know all relative reaction rates of the MM reaction scheme. Additionally, we have the one non-trivial concentration ratio $\beta$ for the given initial substrate and enzyme concentrations.
\end{remark}

\subsubsection*{Stoichiometric compatibility class (conservation laws)}
For the dimensionless system \eqref{eq:dimless-EK}, we identify two linearly independent dimensionless left nullvectors of the stoichiometry matrix $\tilde {\mathcal S}$ given by $l_1=(0,1,1,0)^\top$  and $l_2=(1,0,\beta,1)^\top$ which correspond to the two dimensionless conservation laws\textcolor{black}{\footnote{\textcolor{black}{The original dimensional conservation laws are: $E+C=e_0$ and $S+C+P=s_0$.}}}:
\begin{equation}
    e+c=1\,,\quad s+\beta c +p = 1\,. \label{conservation_laws}
\end{equation}
These two conservation laws allow us to reduce the original four-dimensional problem \eqref{eq:dimless-EK} to a two-dimensional sub-problem restricted to the stoichiometric compatibility class
$$
\mathcal {A}_{\tilde{\mathcal{S}}}=\{(s,e,c,p)\in\mathbb{R}_{\geq 0}^4\,:  e+c=1\,, s+\beta c +p = 1\,\}\,.
$$
We solve these linear conservation laws globally as a graph over, e.g., $(s,c)$-coordinate space:
$$
e=e(s,c)=e(c)=1-c\,,\quad p=p(s,c)=1-s-\beta c\,. 
$$
The corresponding two-dimensional model\footnote{which represents the exact evolution of the original 4D model in the 2D affine subspace/simplex (stoichiometric compatibility class) of the phase space $\mathbb{R}^4_{\geq 0}$. Note also that $\mathcal{V}^\top$ is chosen to correspond to the first and third ODE of system \eqref{eq:dimless-EK} and $\mathcal{U}$ is then calculated accordingly so that $\mathcal{\tilde{S}} = \mathcal{U} \mathcal{V}^\top$.} is given by
\begin{equation}\label{eq:dimless-EK-2D}
\begin{aligned}
\begin{pmatrix}
\frac{ds}{dt}\\

\frac{dc}{dt}
\end{pmatrix}
 &= \mathcal{V}^\top \tilde{V}(y) = 
\begin{pmatrix}
-\beta & +\beta & 0 & 0\\
+1 & -1 & -1 & +1
\end{pmatrix}
\,
\begin{pmatrix}
s(1-c)\\
\alpha c \\
\gamma c\\
\delta (1-s-\beta c)(1-c)
\end{pmatrix}
\\
 &= 
    \begin{pmatrix}
        -\beta \\ +1
    \end{pmatrix} 
    s(1-c)
    + \begin{pmatrix}
        \beta \\ -1 
    \end{pmatrix} \alpha c 
    + \begin{pmatrix}
        0 \\ -1
    \end{pmatrix} \gamma c  + 
    \begin{pmatrix}
    0\\+1
    \end{pmatrix}
    \delta (1-s-\beta c)(1-c)
\end{aligned}
\end{equation}
%
with initial conditions $s(0)=1$ and $c(0)=0$, where we emphasize the splitting into the four reactions involved in the reaction scheme with their reduced dimensionless `stoichiometry' vectors. The irreversible MM reaction scheme ($\delta = 0$) ODEs are

\begin{equation}\label{eq:dimless-EK-2Dirr}
\begin{aligned}
\begin{pmatrix}
\frac{ds}{dt}\\

\frac{dc}{dt}
\end{pmatrix}
 &= 
    \begin{pmatrix}
        -\beta \\ +1
    \end{pmatrix} 
    s(1-c)
    + \begin{pmatrix}
        \beta \\ -1 
    \end{pmatrix} \alpha c 
    + \begin{pmatrix}
        0 \\ -1
    \end{pmatrix} \gamma c 
\end{aligned}
\end{equation}

\subsection{Pre-processing the KF reaction scheme}
The KF model arises from a simplified circadian clock model described by a gene-regulatory network (GRN) with sequestration, see e.g., \cite{kimforger,kimtyson,kim2021}. This GRN can be recast as a CRN (see e.g. \cite{gene3,gene2,gene1}) as follows

\begin{align}
\ce{S ->[k_1] X + S}, \quad &\ce{X ->[k_2] },\quad \ce{X ->[k_3] X + Y} \nonumber \\
&\ce{Y ->[k_4] },\quad \ce{Y ->[k_5] Y + Z} \nonumber \\
&\ce{Z ->[k_6] },\quad \ce{S + Z <=>[k_f][k_b] C}, \quad \ce{C ->[k_u] S}.\label{KFcrn}
\end{align}
We call this CRN as the KF reaction scheme. It describes a transcriptional negative feedback loop for the mRNA $X$, encoded protein $Y$ and repressor protein $Z$. The negative feedback is due to sequestration, where the repressor protein $Z$ suppresses the activator $S$. 

\begin{remark} \label{KF_assumption}
    The KF model is derived under the assumption that $k_u = k_6$, followed by a tQSSA, see e.g. \cite{kim2021}.
    In our analysis, we allow for $k_u \neq k_6$.
\end{remark}
\subsubsection*{The law of mass action}
For CRN \eqref{KFcrn}, the LMA gives the following system of ODEs
\begin{align}
\begin{pmatrix}
\frac{dX}{dT} \\
\frac{dY}{dT} \\
\frac{dZ}{dT} \\
\frac{dS}{dT} \\
\frac{dC}{dT} 
\end{pmatrix} =   \begin{pmatrix}
 1 & -1 & 0 & 0 & 0 & 0 & 0 & 0 & 0\\
 0 & 0 & 1 & -1 & 0 & 0 & 0 & 0 & 0\\
 0 & 0 & 0 & 0 & 1 & -1 & -1 & 1 & 0\\
 0 & 0 & 0 & 0 & 0 & 0 & -1 & 1 & 1\\
 0 & 0 & 0 & 0 & 0 & 0 & 1 & -1 & -1
\end{pmatrix}  
\begin{pmatrix}
k_1 S \\ 
k_2 X  \\
k_3 X  \\
k_4 Y \\
k_5 Y \\
k_6 Z \\
k_f SZ \\
k_b C \\
k_u C
\end{pmatrix} = S V(\hat{X}) \label{KF-full}
\end{align}
with ICs $(X(0),Y(0),Z(0),S(0),C(0)) = (x_0,y_0,z_0,s_0,0)$, where $x_0,y_0,z_0,s_0>0$ and the reaction rates are also strictly positive. Since the sequestration reaction is the first reaction of the MM reaction scheme \eqref{crn2}, the units for $S,Z,T$ and the parameters $k_f$ and $k_b$ are the same as the corresponding variables and parameters in Table \ref{units_table}.

\begin{table}[ht]
\centering
\begin{tabular}{ |c|c|c|c| } 
 \hline
Variable/parameter & Unit & Parameter & Unit \\ \hline
X & M (mol$\cdot$L$^{-1}$) & $k_i$, $i=\{1,\cdots,6\}$ & s$^{-1}$ \\ 
Y & M & $k_u$ & s$^{-1}$ \\ \hline
\end{tabular}
\caption{Units of the variables and parameters for the remaining variables and parameters of the KF reaction scheme.} \label{units_table_KF}
\end{table}

\begin{remark}
    Typically, the first entry of $V(\hat{X})$ is written as $\frac{k_1S}{s_0}$, i.e. the effect on $X$ is proportional to $\frac{S}{s_0}$; see e.g \cite{kimforger,kimtyson,kim2021}. However, here we absorb the $\frac{1}{s_0}$ into $k_1$.
\end{remark}

\subsubsection*{Dimensional analysis}
We define the rescaling of the species concentration as follows
$$
X=k_X x,\quad Y=k_Y y,\quad Z=k_{Z} z,\quad S = k_S s, \quad C = k_C c,\quad T=k_T t.
$$
The reference scales we have chosen are
\begin{align*}
k_S = s_0, \quad k_X = x_0, \quad k_Y = y_0, \quad k_Z = k_C = z_0, \quad k_T = \frac{1}{k_f s_0}.
\end{align*}
The nondimensional parameters are
\begin{align*}
\alpha = \dfrac{k_b}{k_f s_0}, \,\,\, \beta = \dfrac{z_0}{s_0}, \,\,\, \rho_1 = \dfrac{k_1}{k_{f} x_0}, \,\,\, \rho_3 = \dfrac{k_3 x_0}{k_{f}s_0 y_0}, \,\,\, \rho_5 = \dfrac{k_5 y_0}{k_{f}s_0 z_0}, \,\,\, \rho_{2i} = \dfrac{k_{2i}}{k_{f} s_0}, \,\,\, \gamma = \dfrac{k_u}{k_f s_0}
\end{align*}
for $i=1,2,3$, with similar interpretations as in the MM reaction scheme. System \eqref{KF-full} then becomes
\begin{align}
\begin{pmatrix}
\frac{dx}{dt} \\
\frac{dy}{dt} \\
\frac{dz}{dt} \\
\frac{ds}{dt} \\
\frac{dc}{dt} 
\end{pmatrix} =   \begin{pmatrix}
 1 & -1 & 0 & 0 & 0 & 0 & 0 & 0 & 0 \\
 0 & 0 & 1 & -1 & 0 & 0 & 0 & 0 & 0 \\
 0 & 0 & 0 & 0 & 1 & -1 & -1 & 1 & 0 \\
 0 & 0 & 0 & 0 & 0 & 0 & -\beta & \beta & \beta \\
 0 & 0 & 0 & 0 & 0 & 0 & 1 & -1 & -1
\end{pmatrix}  
\begin{pmatrix}
\rho_1 s \\ 
\rho_2 x  \\
\rho_3 x  \\
\rho_4 y \\
\rho_5 y \\
\rho_6 z \\
 sz \\
\alpha c \\
\gamma c
\end{pmatrix} = \tilde {\mathcal S} \tilde V(\hat{x}). \label{KF-full-nondim}
\end{align}
\subsubsection*{Stoichiometric compatibility class (conservation laws)}
We identify one left null-vector of the stoichiometry matrix $\tilde{\mathcal{S}}$ given by $l_1 = (0,0,0,1,\beta)$, corresponding to the dimensionless conservation law
\begin{align}
s + \beta c = 1. \label{conserved-KF}
\end{align}
This allows us to reduce the five-dimensional problem \eqref{KF-full-nondim} to a four-dimensional sub-problem restricted to the stoichiometry compatibility class
\begin{align*}
    \mathcal {A}_{\tilde{\mathcal{S}}} = \{ (x,y,z,s,c) \in \mathbb{R}^5_{\geq 0} \, : \, s + \beta c = 1 \}.
\end{align*}
We can solve the linear conservation law globally as a graph over the $s$-coordinate space,
\begin{align*}
    c = c(s) = \beta^{-1} (1 - s).
\end{align*}
The resulting four-dimensional model is then given by
\begin{align}
\begin{pmatrix}
\frac{dx}{dt} \\
\frac{dy}{dt} \\
\frac{dz}{dt} \\
\frac{ds}{dt}
\end{pmatrix} =  \mathcal{V}^\top \tilde{V}(\hat{y}) = \begin{pmatrix}
 1 & -1 & 0 & 0 & 0 & 0 & 0 & 0 & 0 \\
 0 & 0 & 1 & -1 & 0 & 0 & 0 & 0 & 0 \\
 0 & 0 & 0 & 0 & 1 & -1 & -1 & 1 & 0 \\
 0 & 0 & 0 & 0 & 0 & 0 & -\beta & \beta & \beta
\end{pmatrix}  
\begin{pmatrix}
\rho_1 s \\ 
\rho_2 x  \\
\rho_3 x  \\
\rho_4 y \\
\rho_5 y \\
\rho_6 z \\
 sz \\
 \alpha \beta^{-1} (1-s) \\
\gamma \beta^{-1} (1-s)
\end{pmatrix}. \label{KF-4D-nondim}
\end{align}

\section{Model Reductions of CRNs Based on ci-GSPT}\label{sec:GSPT}
This section establishes the framework for our analysis. First, we make assumptions on the relative magnitudes of the dimensionless parameters, which casts the system as a set of singular perturbation problems. Second, we introduce a coordinate-independent toolbox that leverages the resulting timescale separation in system \eqref{LMA-stoich-red} to perform rigorous model reductions.

\subsection{Asymptotic assumptions on dimensionless parameters}
The pre-processing steps in Section \ref{sec:reaction} for a reaction network model described so far leads to a stoichiometrically reduced $r$-dimensional problem \eqref{LMA-stoich-red}.
Depending on the relative order of magnitudes of the $p$ dimensionless system parameters $\mu\in\mathbb{R}_{\geq 0}^p$, system \eqref{LMA-stoich-red} may consist of processes evolving on two or more distinct timescales.

\begin{assumption}\label{parameterassumption}
Let $0 < \varepsilon \ll 1$ denote a sufficiently small parameter so that the system parameters $\mu_i$, $i=1,\ldots p,$ can be assigned to three different categories:
\begin{itemize}
    \item parameters $\mu_{j}\sim \varepsilon$, $j=1,\ldots p_s$; these are  asymptotically \textbf{small} parameters of $\mathcal{O}(\varepsilon)$;
    \item parameters $\mu_{k}\sim \frac{1}{ \varepsilon}$, $k=1,\ldots p_u$; these are  asymptotically \textbf{large} parameters of order $\mathcal{O}\left(\frac{1}{ \varepsilon}\right)$;
    \item the remaining parameters $\mu_{l}$, $l=1,\ldots p-(p_s+p_u)$; these are considered of order $O(1)$.
\end{itemize}
We relate all asymptotically small parameters $\mu_{j}=\mathcal{O}(\varepsilon)$, $j=1,\ldots p_s$, via $\mu_j = \tilde{\mu}_j \varepsilon $ for some $\tilde{\mu}_j=\mathcal{O}(1)$. Similarly, we relate all asymptotically large parameters $\mu_{k}=\mathcal{O}\left(\frac{1}{ \varepsilon}\right)$, $k=1,\ldots p_u$, via $\mu_k = \tilde{\mu}_k /\varepsilon $ for some $\tilde{\mu}_k=\mathcal{O}(1)$. 
\end{assumption}

Under this specific asymptotic parameter assumption, system \eqref{LMA-stoich-red} defines a perturbation problem
\begin{equation}\label{eq:dimless-EK-general}
    \frac{dy}{dt} =\mathcal{F}(y,\varepsilon)= \sum_{i=0}^{j} \varepsilon^i F_i(y)
\end{equation}
where $j\ge 1$, $0<\varepsilon \ll 1$ is a perturbation parameter, $y \in \mathbb{R}^r$, $r \geq2$ are the dependent variables, $t\in \mathbb{R}$ is the independent variable and $\mathcal F(y,\varepsilon): U\times I \subseteq \mathbb{R}^r \times \mathbb{R} \to \mathbb{R}^r$ is a smooth vector field dependent on a single perturbation parameter $\varepsilon$. It is assumed that $\mathcal F(y,0)=F_0(y) \not\equiv 0$. 
{
\begin{remark} \label{rescale_time}
    Since LMA-based kinetics involves only polynomial terms, the vector field $\mathcal{F}(y,\varepsilon)$ is guaranteed to be analytic. Consequently, upon imposing Assumption \ref{parameterassumption}, $\mathcal{F}$ admits the power series representation given in \eqref{eq:dimless-EK-general}. If parameters of order $\mathcal{O}\left(\frac{1}{\varepsilon}\right)$ appear, the leading-order term $F_0(y)$ may evolve on a fast timescale $\mathcal{O}\left(\frac{1}{\varepsilon^m}\right),\,m\ge 1$. In such cases, one rescales time by $dt=\varepsilon^m d\tilde t$ to recover the perturbation problem \eqref{eq:dimless-EK-general} with respect to $\tilde t$. With a slight abuse of notation, we subsequently drop the tilde $\tilde{}$.
\end{remark}}

\subsection{The tools of ci-GSPT}
\label{sec:ciGSPTtools}
{The following provides a concise summary of the theoretical basis for {calculating model reductions via a timescale disparity, which} integrates \emph{coordinate-independent} (ci)-GSPT with the parametrization method. We emphasize that the underlying theory for both components is well-established; see e.g. \cite{goekewalcher,feliukruffwalcher,wechselberger2020} for ci-GSPT, {\cite{param} for the parametrization method, and \cite{multiple} for the integration of both components.} 
While ci-GSPT provides a geometric framework to identify reduced models without ad-hoc coordinate transformations, the parametrization method enables the explicit computation of both the parametrization of the invariant slow manifold and the resulting slow dynamics, including higher-order terms. We refer the reader to the references cited for further details.}

\subsubsection*{The geometric definition of a singular perturbation problem}
To be able to make further rigorous {\em approximative} model reductions based on timescale separation, system \eqref{eq:dimless-EK-general} has to be identified as a singular perturbation problem.

\begin{definition} \label{def_singular}
System \eqref{eq:dimless-EK-general} is a {\em singular perturbation problem} if the zero level set $S=\{F_0(y)=0\}$ contains a $k$-dimensional, locally connected submanifold $S_0 \subseteq S$ ($1\le k< r$), referred to as the {\em critical manifold}. 
\end{definition}

{
\begin{remark} \label{restrict}
    Assumption \ref{parameterassumption} restricts our analysis to parameter regimes spanning three orders of magnitude: $\mathcal{O}(\varepsilon), \mathcal{O}(1)$ and $\mathcal{O}\left(\frac{1}{\varepsilon}\right)$. While singularly perturbed problems of interest may exist outside these regimes, we emphasize that the tools described in this section apply to any singularly perturbed problem of the form \eqref{eq:dimless-EK-general}, {provided the critical manifold remains normally hyperbolic (see Definition \ref{def_normallyhyperbolic})}.
\end{remark}}

\begin{assumption}
    System \eqref{eq:dimless-EK-general} is a singular perturbation problem with leading-order term factored as $F_0(y)=N_0(y)f_0(y)$ with $r\times (r-k)$ matrix $N_0(y)$ and column vector $f_0(y)\in\mathbb{R}^{r-k}$ such that {$N_0(y)$ has full column rank $r$ for all $y\in S_0$ and $S_0=\{f_0(y)=0\}$}. 
\end{assumption}
\noindent
Thus, the singular perturbation problem \eqref{eq:dimless-EK-general} under study is then given by
\begin{equation}\label{eq:dimless-EK-generalf}
    \frac{dy}{dt} =\mathcal{F}(y,\varepsilon)= N_0(y)f_0(y) + \sum_{i=1}^{j} \varepsilon^i F_i(y) = N_0(y)f_0(y) + G(y,\varepsilon)\,.
\end{equation}

\subsection*{The layer problem}
The singularly perturbed system \eqref{eq:dimless-EK-generalf} is considered the {\em fast(est) system} since it evolves on the fast(est) timescale $t$ identified. Taking the (singular) limit $\varepsilon\to 0$ in \eqref{eq:dimless-EK-generalf} gives the so-called {\em layer problem}
\begin{equation}\label{eq:dimless-general-layer}
    \frac{dy}{dt} =\mathcal{F}(y,0)= N_0(y)f_0(y)\,,
\end{equation}
which describes the leading-order fast(est) evolution of \eqref{eq:dimless-EK-generalf}. Importantly, the critical manifold $S_0=\{f_0(y)=0\}$ forms a $k$-dimensional manifold of equilibria for the layer problem \eqref{eq:dimless-general-layer}. 
{
\begin{definition} \label{def_normallyhyperbolic}
    The critical manifold $S_0$ is called {\em normally hyperbolic} if the linearization of system \eqref{eq:dimless-general-layer} along $S_0$ has $k$ zero eigenvalues and $r-k$ eigenvalues bounded away from the imaginary axis. The $k$ zero eigenvalues, corresponding to the $k$-dimensional tangent space $T_z S_0$ at each $z\in S_0$, are the {\em trivial eigenvalues} and the remaining $r-k$ are the {\em non-trivial eigenvalues}. The critical manifold $S_0$ is {\em attracting} ({\em repelling}) if all
    non-trivial eigenvalues have negative (positive) real parts.
\end{definition}
\begin{lemma} \label{non-trivial_eig}
    The non-trivial eigenvalues of the linearization of system \eqref{eq:dimless-general-layer} along the $k$-dimensional critical manifold $S_0$ are given by the eigenvalues of the $(r-k)\times (r-k)$ matrix $D\!f_0 N_0|_{S_0}$.
\end{lemma}
\begin{proof}
    See, e.g., the proof of Lemma 3.3 in \cite{wechselberger2020}.
\end{proof}}

\begin{definition}
    In system \eqref{eq:dimless-general-layer},
    the column vectors of $N_0(y)$ attached to their base point $y\in S_0$ form a {\em fast (linear) fiber bundle} $\mathcal{N}_0=\cup_{y\in S_0} N_0(y)$ of the critical manifold $S_0$, i.e., they form the corresponding non-trivial eigenbasis to the non-trivial eigenvalues.
\end{definition}
\begin{remark}
    Since the critical manifold $S_0$ is a manifold of equilibria for the layer problem with $\dim S_0 = k$, there are always $k$ trivial zero eigenvalues along $S_0$. The union of the corresponding nullspaces attached to the base points forms the tangent bundle $TS_0=\cup _{y\in S_0} T_y S_0$.
\end{remark}

In the normally hyperbolic attracting case, the layer problem \eqref{eq:dimless-general-layer} describes the (non-linear) evolution towards the critical manifold $S_0$, 
i.e., for a given initial condition $y(0)=y_0$ near $S_0$, the initial value problem evolves (exponentially fast) towards its base point $y_0^b \in S_0$.
At each base point $y_0^b \in S_0$, this $(r-k)$-dimensional non-linear layer flow is tangent to $N_0(y)$, i.e., the fast fiber bundle $\mathcal{N}_0$ represents the linear approximation of fast (non-linear) fiber bundle $\mathcal{W}^s(S_0)$. 

Hence, the (trivial) tangent bundle $T\mathbb{R}^r$ of the phase space $\mathbb{R}^r$ restricted to the critical manifold $S_0$ has the splitting 
    $$T\mathbb{R}^r = TS_0 \oplus \mathcal N_0$$ 
    which encodes the dynamic information (slow and fast directions) of the singular perturbation problem near $S_0$; see Figure \ref{fig:param}.

\subsubsection*{The reduced problem}

\begin{figure}[ht]
\centering
\includegraphics[scale=0.5]{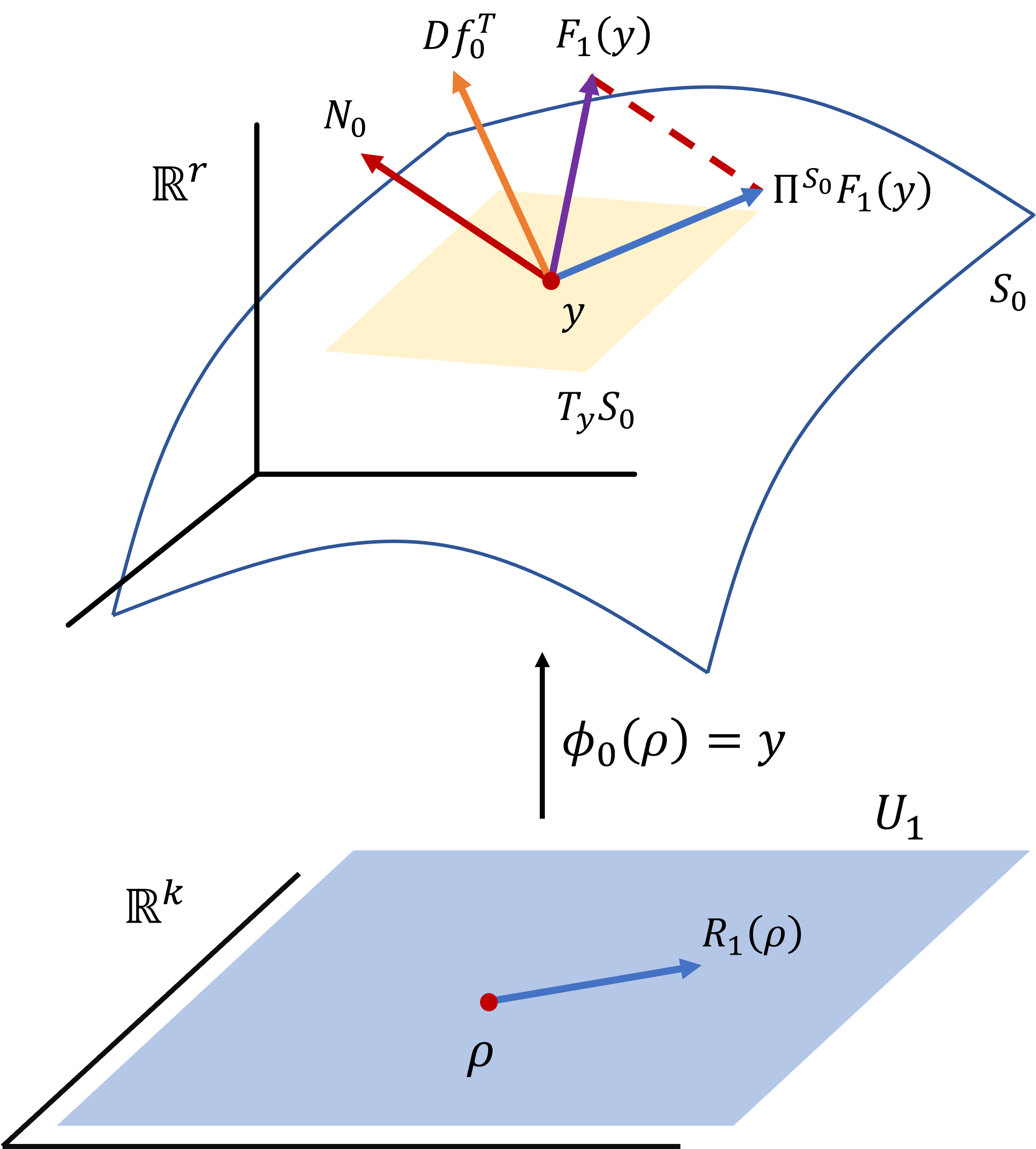}
\caption{A sketch of the {graph} embedding $\phi_0$ of the critical manifold $S_0$ and  the projection of the vector field $F_1(y)$ from $T_y\mathbb{R}^r$ onto $T_yS_0$.}
\label{fig:param}
\end{figure}

In order to find the leading-order approximation of the slow dynamics near $S_0$, we begin by rescaling time $\tau = \varepsilon t$ in \eqref{eq:dimless-EK-generalf} to obtain the \textit{slow system}

\begin{equation}\label{eq:dimless-EK-generals}
    \frac{dy}{d\tau} = \dfrac{1}{\varepsilon} N_0(y) f_0(y) + \sum_{i=0}^{j-1} \varepsilon^i F_{i+1}(y)
\end{equation}
evolving on the slow timescale $\tau$. We emphasize that {for $\varepsilon > 0$} the fast and the slow system are equivalent, i.e., they describe the same problem on different timescales. The singular nature of the slow problem becomes apparent when taking the limit $\varepsilon \to 0$ in system \eqref{eq:dimless-EK-generals} which is only well-defined if
\begin{itemize}
    \item the variable $y$ is constrained to the critical manifold, i.e., $y\in S_0$, and
    \item the vector field $F_1$ along the critical manifold $S_0$ is constrained to its slow component in the tangent bundle $TS_0$.
\end{itemize}
This desired slow component of $F_1$ can be obtained via a uniquely defined oblique projection $\Pi^{S_0}: T \mathbb{R}^r \to T S_0$ along the fast linear fiber bundle $\mathcal{N}_0$ onto $TS_0$; see Figure~\ref{fig:param}. Taking the limit in this careful way results in the desired reduced problem on $S_0$ which is a differential-algebraic problem 

\begin{equation} \label{reduced:gen}
\dfrac{dy}{d\tau} = \Pi^{S_0} F_1(y)=\left(\mathbb{I}_r - N_0(y) \left(D\!f_0(y) N_0(y)\right)^{-1} D\!f_0(y)\right)F_1(y) ,\qquad f_0(y) =0,
\end{equation}
where $\mathbb{I}_r$ denotes the $r \times r$ identity matrix. Analysis of such a differential-algebraic problem  is done in local coordinate charts of the critical manifold $S_0$.

{
\begin{assumption} \label{graph}
     The critical manifold $S_0$ is given by a graph embedding $\phi_0 : U_1 \subseteq \mathbb{R}^k  \to \mathbb{R}^r$, defined by $\phi_0(\rho) = (\rho, \psi_0(\rho))^\top$, where the coordinates $y \in \mathbb{R}^r$ are {partitioned} as $y= (\rho, \eta)^\top$ with $\rho \in \mathbb{R}^k$; see Figure \ref{fig:param}.
\end{assumption}}
In the coordinate chart $\rho\in U_1$, we seek the vector field $R_1(\rho)$ whose image under the linear mapping $D\phi_0$ coincides with the slow component vector field $\Pi^{S_0} F_1 \in TS_0$, i.e.,
$$D\phi_0 R_1(\rho) = \Pi^{S_0} F_1(\phi_0(\rho))\,,\quad \forall \rho\in U_1\,.$$ 
This yields the leading-order slow flow approximation, defining the $k$-dimensional {\em reduced problem} in the coordinate chart $\rho\in U_1 \subseteq \mathbb{R}^k$:
\begin{equation}\label{reduced:gen-param}
   \dfrac{d\rho}{d\tau}= R_1(\rho):=D\phi_0^L \Pi^{S_0} F_1(\phi_0(\rho))\,,
\end{equation}
where $D\phi_0^L$ denotes a left inverse of the $r\times k$ matrix $D\phi_0$.  
{Given the graph structure of $\phi_0$, {we set $D\phi_0^L = \left(\mathbb{I}_k \,\, \mathbb{O}_{k,r-k} \right)$\footnote{{Note, left inverses are non-unique.}} to achieve a projection of $\Pi^{S_0} F_1$ onto the coordinate chart $\rho\in U_1$.}

\subsubsection*{Tikhonov-Fenichel theory and the parametrization method}

Tikhonov-Fenichel theory {provides a rigorous mathematical foundation} for analyzing systems with multiple timescales:
\begin{theorem}[Persistence of the slow manifold \cite{fenichel}; see also \cite{tikhonov}]
    Given a singularly perturbed system \eqref{eq:dimless-EK-generalf} with an analytic vector field and a compact normally hyperbolic critical manifold $S_0$ (possibly with boundary). {Then, for any finite integer $k \ge 1$, there exists an $\varepsilon_0 = \varepsilon_0(k) > 0$ such that for all $0 < \varepsilon \le \varepsilon_0$, the following hold:}
    \begin{enumerate}
        \item There exists a {$C^k$} differentiable manifold $S_0^\varepsilon$, locally invariant under the flow of \eqref{eq:dimless-EK-generalf}, which is $\mathcal{O}(\varepsilon)$ close to $S_0$ and possesses a unique asymptotic expansion in $\varepsilon$ (up to order $k$).
        \item The manifold $S_0^\varepsilon$ is generally not unique, but all such manifolds differ by terms of order $\mathcal{O}(e^{-c/\varepsilon})$ for some $c > 0$.
        \item The corresponding slow flow on the manifold $S_0^\varepsilon$ is of class {$C^k$} and converges to the reduced flow on the critical manifold $S_0$ as $\varepsilon \to 0$.
    \end{enumerate}
\label{fenichel_theorem}
\end{theorem}

\begin{remark}\label{rem:uniqueness}
{
    While the vector field in \eqref{eq:dimless-EK-generalf} is analytic, the slow manifold $S_0^\varepsilon$ guaranteed by Theorem \ref{fenichel_theorem} is generally not analytic. Rather, for any finite $k$, one can find a $C^k$ manifold, but it is typically not unique. Specifically, in the analytic case, different realizations of $S_0^\varepsilon$ differ only by exponentially small terms of $\mathcal{O}(e^{-c/\varepsilon})$ for some $c>0$.} 

    {
    Crucially, all such candidate manifolds share the \emph{same} unique asymptotic power series expansion in $\varepsilon$. The parametrization method computes the coefficients of this unique power series order-by-order. Consequently, for the finite-order approximations typically used in applications, the reduced model we derive is the valid approximation for \emph{any} realization of the slow manifold $S_0^\varepsilon$, as the non-uniqueness lies beyond all orders of the asymptotic expansion.
 }   
\end{remark}

\begin{remark} \label{eq_remark}
    A hyperbolic equilibrium $\rho^*$ of the reduced problem \eqref{reduced:gen} persists under the perturbation of $0<\varepsilon\ll 1$. It is located on the slow manifold $S_0^\varepsilon$ and by Theorem \ref{fenichel_theorem}, is a hyperbolic equilibrium of the full system \eqref{eq:dimless-EK-generalf}.
\end{remark}

{
We now introduce the \emph{parametrization method} as the constructive counterpart to Tikhonov-Fenichel theory. While the theory guarantees the existence of the slow manifold and its unique asymptotic expansion, the parametrization method provides the {framework} to explicitly calculate these series coefficients to any order. This general {method} is well-established; see e.g. the monograph \cite{param}. In the following, we adopt the specific implementation for multiple timescale systems developed in \cite{multiple}.}\\

Based on Theorem \ref{fenichel_theorem}, $S_0$ perturbs to nearby $S_0^\varepsilon$ which is (locally) invariant under the flow of the full system \eqref{eq:dimless-EK-generalf}. 
Since $S_0$ is a graph $(\rho,\psi_0(\rho))$, we may search for a graph embedding $\phi:  \mathbb{R}^k \times I \to \mathbb{R}^r$, $I := [0,\varepsilon_0)$, of $S_0^\varepsilon$, with the power series expansion
\begin{align} \label{SM_graph}
    \phi(\rho,\varepsilon) = \begin{pmatrix}
        \rho \\ \psi_0(\rho)
    \end{pmatrix} + \varepsilon \begin{pmatrix}
        0 \\ \psi_1(\rho)
    \end{pmatrix} +  \varepsilon^2 \begin{pmatrix}
        0 \\ \psi_2(\rho)
    \end{pmatrix} + \cdots = \begin{pmatrix}
        \rho \\ \psi_0(\rho)
    \end{pmatrix} + \begin{pmatrix}
        0 \\ \sum_{j=1} 
    \varepsilon^j \psi_j(\rho)
    \end{pmatrix}.
\end{align}
\noindent
Theorem \ref{fenichel_theorem} also states that the corresponding slow flow on $S_0^\varepsilon$ converges to the reduced flow on $S_0$. Thus, we make a power series ansatz for the slow vector field
\begin{align}
    \dfrac{d\rho}{dt} = \varepsilon R(\rho,\varepsilon) = \varepsilon R_1(\rho) + \varepsilon^2 R_2(\rho) + \varepsilon^3 R_3(\rho) + \cdots =\sum_{j=1} \varepsilon^j R_j(\rho) \label{slow_vector_field}
\end{align}
in the coordinate chart $\rho\in\mathbb{R}^k$, where we have written system \eqref{slow_vector_field} in the original fast timescale $t=\tau/\varepsilon$.
\begin{remark} \label{remark_ICs}
    To obtain a well-posed initial value problem on $S_0^\varepsilon$ with the vector field \eqref{slow_vector_field}, one needs to determine the corresponding base point on $y_{0,\varepsilon}^b\in S_0^\varepsilon$ of the fast fiber containing the given initial condition $y_0$. By Tikhonov-Fenichel theory \cite{tikhonov,fenichel}, the fast fiber bundle $\mathcal{W}^{s}(S_0)$ of the layer problem also perturbs to a nearby fast fiber bundle $\mathcal{W}^{s}(S_0^\varepsilon)$ and the foliation of $\mathcal{W}^{s}(S_0^\varepsilon)$  is invariant in the sense that $\mathcal{W}^{s}(y_{0,\varepsilon}^b) \cdot t \subset \mathcal{W}^{s}(y_{0,\varepsilon}^b\cdot t)$, such that $y_{0,\varepsilon}^b,y_{0,\varepsilon}^b\cdot t \in S_0^\varepsilon$. This then justifies using the base point as the IC for system \eqref{slow_vector_field}.
\end{remark}

{The invariant slow manifold $S_0^\varepsilon$ satisfies the {\em conjugacy equation}}
\begin{align}
   D\phi(\rho,\varepsilon) R(\rho,\varepsilon)= F(\phi(\rho,\varepsilon),\varepsilon)\,.  \label{conjugacy_equation}
\end{align}
which is, in general, a non-linear partial differential equation (PDE). 
{A power series expansion of the conjugacy equation \eqref{conjugacy_equation} in $\varepsilon$ allows us to calculate an approximate solution by solving iteratively for the expansion coefficients.} The $\mathcal{O}(\varepsilon^j)$ term is given by
\begin{align*}
    D\phi_0(\rho) \cdot R_j(\rho) - N_0(\phi_0(\rho)) Df_0(\phi_0(\rho)) \cdot \phi_j =  G_j(\phi_0, \cdots, \phi_{j-1},r_1,\cdots,r_{j-1})
\end{align*}
and we call this the {\em infinitesimal conjugacy} or {\em homological equation}. The inhomogeneous term $G_j$ involves terms calculated in previous iterations and the first three are given by
\begin{align} \label{G}
       G_1 &= F_1(\phi_0), \nonumber \\
       G_2 &= -D \phi_{1} \cdot R_{1}+F_{2}\left(\phi_{0}\right) +D F_{1}\left(\phi_{0}\right) \cdot \phi_{1}+\frac{1}{2} D^{2} F_{0}\left(\phi_{0}\right)\left(\phi_{1}, \phi_{1}\right), \nonumber \\ 
       G_3 &= -D\phi_2 \cdot R_1 - D \phi_1 \cdot R_2 + \dfrac{1}{6} D^3 F_0 (\phi_0)(\phi_1,\phi_1,\phi_1) + D^2F_0(\phi_0)(\phi_1,\phi_2) \nonumber \\
       &+ \dfrac{1}{2}D^2F_1(\phi_0)(\phi_1,\phi_1)
       + DF_1(\phi_0) \cdot \phi_2 + DF_2(\phi_0) \cdot \phi_1 + F_3 (\phi_0).
\end{align}
Finally, the solution to the homological equation is obtained by applying the projectors $\Pi^{S_0}$, respectively $\Pi^{N_0} :=  \mathbb{I}_r - \Pi^{S_0}$, which separates and solves for the vector field update $R_j$ from the slow manifold update $\psi_j$:
\begin{align}
    R_j(\rho) = D\phi_0^L \Pi^{S_0} G_j, \quad 
    \psi_j(\rho) = - D\!f_0^R D\phi_0^L \Pi^{N_0} G_j \label{parametrization_formula}
\end{align}
where $y = (\rho,\eta) \in \mathbb{R}^r$, {$D\phi_0^L$ a left inverse} and $D\!f_0^R$ a right inverse. We note that $\psi_j(\rho)$ is a correction term for the graph embedding $\phi(\rho,\varepsilon)$ in the direction of $\eta \in \mathbb{R}^{r-k}$, i.e. $\phi_j(\rho) = \left(0 \,\, \psi_j(\rho) \right)^\top$.

In many applications, calculating the leading-order terms $\phi_0$ and $R_1$ are sufficient. However, it might be the case that, e.g., $R_1(\rho) \equiv 0$. This can happen because (i) $F_1$ aligns with the fast fibers $N_0$ along $S_0$ or, trivially, (ii) $F_1\equiv 0$. In these cases, the leading-order approximation of the slow flow is $\mathcal{O}(\varepsilon^2)$ or slower. 

\begin{remark} \label{multiple_ts}
In general, the parametrization method is needed to resolve multiple timescales, i.e., more than two distinct timescales. This becomes particularly important when the slow flow
$$
\frac{d\rho}{dt}= \varepsilon R(\rho,\varepsilon)
$$
on $S_0^\varepsilon$ is also singularly perturbed, i.e., the level set $\{R(\rho,0)=R_1(\rho)=0\}$ contains a $k_1$-dimensional locally connected submanifold, $1\le k_1<k$, along which the corresponding flow is infra-slow. The parametrization method then allows to approximate this infra-slow flow. For further details, we refer the reader to \cite{multiple}.
\end{remark}

\subsection{ci-GSPT model reduction of the MM reaction model \eqref{eq:dimless-EK-2D}}
\label{sec:MM_GSPT}
The applicability of {the} framework relies solely on the existence of a normally hyperbolic critical manifold. Provided this condition is met, the {approach} detailed in Section \ref{sec:ciGSPTtools} yields {an approximate model reduction underpinned by Tikhonov-Fenichel theory and eliminates the ambiguity associated with selecting an appropriate QSSA \textit{a priori}.} In this section, we apply this computational framework to the MM reaction model \eqref{eq:dimless-EK-2D} and illustrate the results with two examples.

\subsubsection*{Asymptotic parameter assumptions}
Based on Assumption \ref{parameterassumption},
there are $3^4=81$ possible asymptotic parameter configurations for system \eqref{eq:dimless-EK-2D}, the non-dimensionalized MM reaction scheme ($3^3=27$ for the irreversible reaction scheme $\delta=0$). $67$ of these give singular perturbation problems (23 for the irreversible reaction scheme $\delta=0$).

\subsubsection*{Normally hyperbolic critical manifolds}
There are $47$ parameter configurations for the reversible reaction scheme (16 for the irreversible reaction scheme $\delta = 0$) where the critical manifold(s) is(are) normally hyperbolic. The critical manifolds for these $47$ normally hyperbolic cases may be written as a graph over $s$, i.e., $c = c_0(s)$ where $\phi_0 = (s ,\, c_0(s))^\top$:
\begin{remark}
Some important observations:
\begin{itemize}
    \item 2 of 47 reversible cases consists of a single critical manifold which is repelling: $\beta = \mathcal{O}(\varepsilon^{-1}), \alpha = \mathcal{O}(\varepsilon), \delta = \mathcal{O}(1)$ and either $\gamma = \mathcal{O}(\varepsilon)$ or $\mathcal{O}(1)$. 
    \item The remaining 45 of the 47 reversible cases have either: (a) a unique critical manifold that is attracting or (b) one of the two critical manifolds is attracting. Section \ref{MM_class} provides a further discussion on the geometries of the critical manifolds.
    \item For cases with a loss of normal hyperbolicity on a critical manifold, the theory we have presented is still applicable to the normally hyperbolic subsets.
\end{itemize}
\end{remark}

\subsubsection*{Model reductions}
Next, we reduce all these 47 reversible (16 irreversible) cases to a 1D-problem in the $s$-coordinate chart, i.e., the ci-GSPT reduction step introduced before gives
\begin{align} \label{reduce}
 \dfrac{ds}{dt} &=  - \varepsilon R(s,\varepsilon)
 \end{align}
which describes the flow on the one-dimensional invariant slow manifold $S_0^\varepsilon$. The negative sign in system \eqref{reduce} reflects the desire that the substrate must deplete in order to form the product in \eqref{crn2} and \eqref{crn}.

Consider then the second conserved quantity in \eqref{conservation_laws}, $s+\beta c +p=1$. Taking the time derivative and restricting the conserved quantity to the slow manifold $c(s,
\varepsilon)=c_0(s) + \varepsilon c_1(s) + \cdots$ leads to 
\begin{align}
 \dfrac{ds}{dt} &=  -\varepsilon R(s,\varepsilon) \nonumber \\
    \dfrac{dp}{dt} &= + \varepsilon \left(\beta D_s c (s,\varepsilon) +1 \right)  R(s,\varepsilon)
    \label{relation_full}
\end{align}
with $D_s c(s,\varepsilon) = \sum_{i=0} \varepsilon^i c_i'(s)$.
{
\begin{remark} \label{not_approx} 
System \eqref{relation_full} describes the invariant flow on the slow manifold $S_0^\varepsilon$. We emphasize that this is not a heuristic approximation, but rather the formal description of the dynamics. As discussed in Remark \ref{rem:uniqueness}, while the manifold itself is defined only up to exponentially small terms, system \eqref{relation_full} represents the unique asymptotic expansion of the flow common to all realizations. 
\end{remark}}

Leading-order approximations via the parametrization method can then be calculated for the substrate depletion and product formation rates: 
{\begin{align}
    \dfrac{ds}{dt} &= - \varepsilon^j R_j(s) \nonumber \\
    \dfrac{dp}{dt} &= + \varepsilon^j (\beta c_0'(s) +1) R_j(s) \label{relation}
\end{align}}
where $R_j(s)$ is the first nonzero term of the slow vector field \eqref{slow_vector_field}. 

\begin{remark} \label{one_step_remark}
     System \eqref{relation_full} may be formally interpreted as a one-step reaction network
\begin{center}
\begin{tabular}{c c c l}
 \multirow{2}{*}{$\overset{l_+(s,\varepsilon)}{\lcirclearrowright} s \overset{k_+(s,\varepsilon)}{\longrightarrow} p$}  
 & & & $\dfrac{ds}{dt} = -k_+(s,\varepsilon) s + l_+(s,\varepsilon) s$ \\ [2.5ex]
 & & & $\dfrac{dp}{dt} = k_+(s,\varepsilon) s$ 
\end{tabular}
\end{center}
with
\begin{align}
    k_+(s,\varepsilon) =  \varepsilon \left(\beta D_s c (s,\varepsilon) +1\right)  
 \frac{R(s,\varepsilon)}{s}, \quad 
 l_+(s) = \varepsilon \beta D_s c (s,\varepsilon)\frac{R(s,\varepsilon)}{s}. 
 \end{align}
Consequently, system \eqref{relation} may also be formally interpreted as a one-step reaction network, where we have obtained leading-order approximations of $k_+(s,\varepsilon)$ and $l_+(s,\varepsilon)$ as follows
 \begin{align}
    k_+(s) =  \varepsilon^j \left(\beta c_0' (s) +1\right) 
 \frac{R_j(s)}{s}, \quad 
 l_+(s) = \varepsilon^j \beta c_0'(s) \frac{R_j(s)}{s}.
 \end{align}
 \end{remark}

\subsubsection*{Examples of model reductions}
\paragraph{The classical tQSSA: $\varepsilon := \gamma$, $ \alpha,\beta = \mathcal{O}(1)$ and $\delta = 0$} 
The model reduction for $\gamma \ll 1$ has historically been dealt with a tQSSA, i.e., by first introducing the total substrate concentration {$\bar{s}:=s + \beta c$} and then applying a QSSA for $c$; see e.g. \cite{borghans,tzafriri}. In contrast to the typical approach, we use our ci-GSPT tools in order to do a coordinate-independent reduction. 

{ 
\begin{remark}
This tQSSA assumption also fits the criteria for a \textit{partial equilibrium approximation (PEA)}, where one of the reversible reactions is assumed to be fast \cite{schauerandheinrich,goussis}.    
\end{remark}}

System \eqref{eq:dimless-EK-2Dirr} for the given parameter configuration is
\begin{align*}
\begin{pmatrix}
\frac{ds}{dt} \\ 
\frac{dc}{dt}
\end{pmatrix} =N_0 f_0(s,c)  + \varepsilon F_1(s,c) 
= \begin{pmatrix}
-\beta \\ 1
\end{pmatrix} (s(1-c)-\alpha c) + \varepsilon \begin{pmatrix} 0  \\
-  c
\end{pmatrix}.
\end{align*}
 The set {$S_0 = \{(s,c) \, | \, f_0(s,c) = 0 \}$} is a one-dimensional manifold and is a graph over $s$:
 \begin{align*}
c_0(s) = \dfrac{s}{\alpha +  s}
\end{align*} 
and, hence, we can write a graph embedding {$\phi(s,\varepsilon)$} of the slow manifold $S_0^{\varepsilon}$:
\begin{align*}
\phi(s,\varepsilon) = \phi_0(s) +  \cdots = \begin{pmatrix}
 s \\ \frac{s}{\alpha + s}
\end{pmatrix}  + \cdots.
\end{align*}
The non-trivial eigenvalue is given by 
{\begin{align*}
Df_0 N_0 \big|_{S_0} =  \begin{pmatrix}
1 - c & -\alpha - s
\end{pmatrix} \begin{pmatrix}
-\beta \\ 1
\end{pmatrix} \bigg|_{S_0} = -\beta\left(1-\frac{s}{\alpha + s}\right) - \alpha - s.
\end{align*}}
We then have the following: 
\begin{align*}
{\Pi}_0^S =  \begin{pmatrix}
\frac{(\alpha+s)^2}{\alpha \beta + (\alpha + s)^2} & \frac{\beta(\alpha+s)^2}{\alpha \beta + (\alpha + s)^2} \\
\frac{\alpha}{\alpha \beta + (\alpha + s)^2} & \frac{\alpha \beta}{\alpha \beta + (\alpha + s)^2} 
\end{pmatrix}\,,\quad 
D\phi_0^L  = \begin{pmatrix}
1 & 0 
\end{pmatrix}\,.
\end{align*}
Applying \eqref{reduced:gen-param} with $G_1 = F_1$, which is schematically shown in Figure \ref{T1_diagram}a, 
we obtain the reduced vector field on $S_0$:
\begin{align}
\dfrac{ds}{dt} =  {\varepsilon} D\phi_0^L {\Pi}_0^S G_1 = -\varepsilon R_1(s) = - \varepsilon \dfrac{ \beta s ( \alpha  +s)}{\alpha \beta + (\alpha + s)^2}\,. \label{tQSSA_model}
\end{align}
Numerical simulations in Figure \ref{T1_diagram}b show good agreement between the 2D system \eqref{eq:dimless-EK-2Dirr} and the model reduction \eqref{tQSSA_model}. Finally, we can also calculate {the} leading-order product formation
{
\begin{align}
\dfrac{dp}{dt} = \varepsilon (\beta c_0'(s) + 1) R_1(s) =  \varepsilon \dfrac{s}{\alpha +  s} \label{product_tqssa}.
\end{align}}

\begin{remark}
    We highlight the asymmetric form for the substrate depletion dynamics and the {product} formation dynamics (i.e. $\frac{ds}{dt} \neq - \frac{dp}{dt}$) which is a hallmark of the tQSSA application in the original $(s,c)$-coordinate system.
\end{remark}

\begin{remark}
From the last equation in system \eqref{eq:dimless-EK}, the leading-order of the dimensionless product formation rate is also given by
\begin{align*}
    \dfrac{dp}{dt} = \varepsilon \beta c_0(s)= \varepsilon \beta \dfrac{s}{\alpha + s} 
\end{align*}  
which matches \eqref{product_tqssa}.
\end{remark}

\begin{figure}[ht]
\centering
\begin{subfigure}{0.45\textwidth}
\centering
 \includegraphics[scale=0.75]{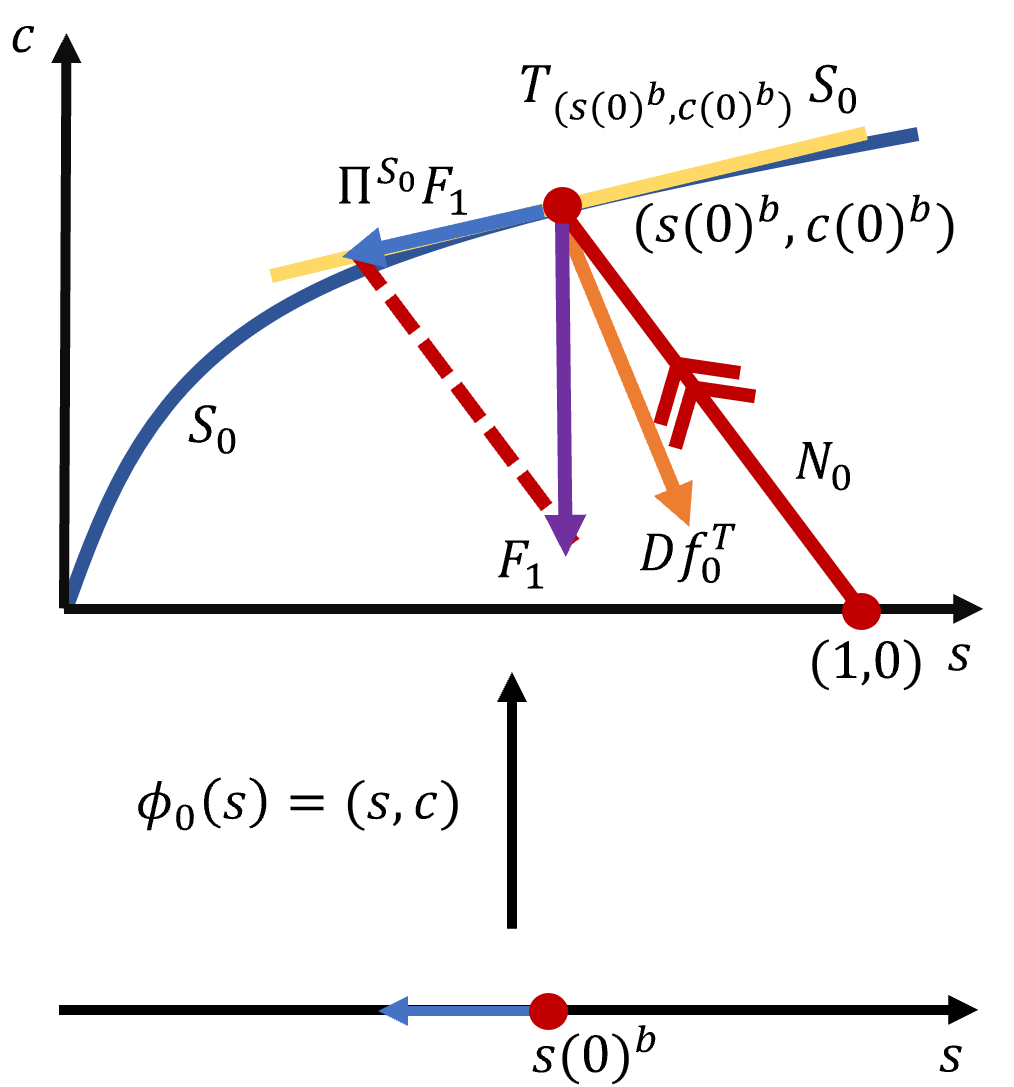}
 \subcaption{}
\end{subfigure}
\begin{subfigure}{0.52\textwidth}
\centering
  \includegraphics[width=1\linewidth]{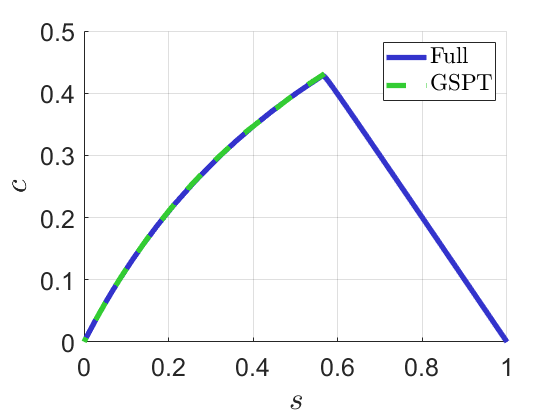}
  \subcaption{}
\end{subfigure}
\caption{(a) A diagram of the ci-GSPT based model reduction of system \eqref{eq:dimless-EK-2Dirr} for $\gamma = \mathcal{O}(\varepsilon), \alpha,\beta = \mathcal{O}(1)$. (b) A comparison of the simulation of the full system \eqref{eq:dimless-EK-2Dirr} and the tQSSA \eqref{tQSSA_model} on its critical manifold. Parameter values are $\alpha = 0.75, \beta = 1, \gamma = 0.005$ and the IC for the model reduction is $s(0) = 0.5687$.} 
\label{T1_diagram}
\end{figure}

\begin{remark} \label{remark_IC_tQSSA}
Since $N_0$ is a constant vector and $f_0$ is a scalar, then the fast linear fiber bundle $\mathcal{N}_0$ and the fast non-linear fiber bundle $\mathcal{W}^S$ coincide. Therefore, the IC of the model reduction is obtained by the intersection of the {layer fibre $N_0$ through $(1,0)$},
\begin{align*}
    c = \beta^{-1}(1-s)\,,
\end{align*}
and the critical manifold $c_0(s)$. Hence, the IC is given by
\begin{align*}
    s(0)^b = \dfrac{-(\alpha + \beta - 1) + \sqrt{(\alpha + \beta - 1)^2 + 4 \alpha }}{2}, \quad c(0)^b = \dfrac{s(0)^b}{s(0)^b+\alpha}.
\end{align*}
The IC $s(0)^b$ is used in the simulation of system \eqref{tQSSA_model} in Figure \ref{T1_diagram}b.
\end{remark}

\paragraph*{An sQSSA case with multiple timescales: $\varepsilon := \beta$, $\alpha = \mathcal{O}(1)$, 
{$\gamma =\varepsilon \tilde\gamma$ with $\tilde\gamma=O(1)$,} and $\delta = 0$}
We now show a case where the calculation of $R_2$ in the slow vector field \eqref{slow_vector_field} is required, demonstrating the use of the parametrization method. 

System \eqref{eq:dimless-EK-2Dirr} with the above parameter configuration is as follows:
\begin{align*}
\begin{pmatrix}
\frac{ds}{dt} \\ 
\frac{dc}{dt}
\end{pmatrix} = N_0 f_0(s,c) + \varepsilon F_1(s,c) =
\begin{pmatrix}
0 \\ 1
\end{pmatrix} (s(1-c)-\alpha c) + \varepsilon \begin{pmatrix} \alpha c -s(1-c)  \\
-  \tilde{\gamma} c
\end{pmatrix}.
\end{align*}
One shows that $R_1 \equiv 0$, i.e., there is no $\mathcal{O}(\varepsilon)$ contribution to the slow flow.\footnote{This is left as an exercise for the reader.} Now, in order to calculate the $O(\varepsilon^2)$ contribution $R_2$, we first require the $\phi_1(s)$ update to the graph embedding $\phi(s,\varepsilon)$,
{
\begin{align}
  c_1(s) = -\dfrac{\tilde{\gamma} s}{(\alpha+s)^2}\,,
\end{align}
which leads to an improved approximation for $S_0^\varepsilon$:
\begin{align}
    \phi(s,\varepsilon) = \begin{pmatrix}
        s \\ \frac{s}{\alpha + s}
    \end{pmatrix} + \varepsilon \begin{pmatrix}
        0 \\ -\frac{\tilde{\gamma} s}{(\alpha+s)^2}
    \end{pmatrix} + \cdots.
\end{align}}
{We then calculate:
\begin{align*}
    \dfrac{ds}{dt} &= \varepsilon^2 D\phi_0^L \Pi^{S_0} G_2(R_1, \phi_0,\phi_1) 
    = \varepsilon^2\begin{pmatrix}
        1 & 0
    \end{pmatrix}\begin{pmatrix}
        1 & 0 \\
        \frac{\alpha}{(\alpha+s)^2} & 0
    \end{pmatrix} G_2(R_1, \phi_0,\phi_1) \nonumber \\
    &=-\varepsilon^2 R_2(s) = -\varepsilon^2\dfrac{\tilde{\gamma} s}{\alpha+ s} 
\end{align*}}
where the four terms of {the} inhomogeneity $G_2$ are:
\begin{align*}
    G_2 &= - D \phi_1 \cdot r_1 + F_2(\phi_0) + D F_1(\phi_0) + \frac{1}{2} D^2 F_0(\phi_0)(\phi_1,\phi_1)  \\
    &= \begin{pmatrix}
        0 \\ 0
    \end{pmatrix} +  \begin{pmatrix}
        0 \\ 0
    \end{pmatrix} + \begin{pmatrix}
        -\frac{\tilde{\gamma} s}{\alpha + s} \\ \frac{\tilde{\gamma}^2 s}{(\alpha + s)^2}
    \end{pmatrix} + \begin{pmatrix}
        0 \\ 0
    \end{pmatrix}.
\end{align*}
The leading-order product formation rate is then given by
{
\begin{align}
\dfrac{dp}{dt} = \varepsilon^2 R_2(s) =  \varepsilon^2 \dfrac{\tilde{\gamma} s}{\alpha +  s}.
\end{align}}

\paragraph*{More examples} Three more `non-standard' model reduction (i.e. using ci-GSPT tools) examples are found in the {Supplementary Material II}: (a) the classical sQSSA $\varepsilon:= \beta, \alpha,\gamma = \mathcal{O}(1), \delta = 0$, (b) the classical rQSSA $\varepsilon:= \beta^{-1}, \alpha,\gamma = \mathcal{O}(1), \delta = 0$ and (c) a reversible sQSSA $\varepsilon:= \beta, \alpha,\gamma,\delta = \mathcal{O}(1)$.

\subsection{ci-GSPT model reduction of the KF reaction model \eqref{KF-full-nondim}}
\label{KF}
In contrast to considering all 81 parameter configurations for the MM reaction scheme, we consider only one possible parameter configuration for the KF reaction scheme to illustrate a model reduction of a larger CRN.

\paragraph*{KF reaction scheme tQSSA: $\varepsilon := \rho_1$, {$\gamma = \varepsilon \tilde{\gamma}, \rho_{i} = \varepsilon \tilde{\rho}_i$ for $i = \{2,3,4,5,6\}$ and $\alpha,\beta = \mathcal{O}(1)$}} Similar to MM, the classical approach is to apply the coordinate transformation $\bar{z} = z + c$ and then to apply the QSSA for $s$; see e.g. \cite{kimtyson,kim2021}. Using ci-GSPT, we avoid a coordinate transformation.

System \eqref{KF-4D-nondim} for this parameter configuration becomes
\begin{align}
\begin{pmatrix}
\frac{dx}{dt} \\
\frac{dy}{dt} \\
\frac{dz}{dt} \\
\frac{ds}{dt} 
\end{pmatrix} =
\begin{pmatrix}
0 \\ 0 \\ \beta^{-1} \\ 1
\end{pmatrix} (- \beta sz + \alpha (1-s))+ \varepsilon\begin{pmatrix}
 s - \tilde{\rho}_2 x  \\ \tilde{\rho}_3 x - \tilde{\rho}_4 y \\ \tilde{\rho}_5 y - \tilde{\rho}_6 z  \\ \tilde{\gamma} (1-s)
\end{pmatrix}.
\end{align}
The set {$S_0 = \{(x,y,z,s) \, | \, f_0(x,y,z,s) = 0 \}$} is a three-dimensional manifold and is a graph over $(x,y,z)$:

\begin{align*}
s(x,y,z) =  s(z) = \dfrac{\alpha}{\alpha + \beta z}.
\end{align*}
The single non-trivial eigenvalue is given by 
{
\begin{align*}
Df_0 N_0 \big|_{S_0} =  \begin{pmatrix}
0 & 0 & -\beta s & -\alpha - \beta z
\end{pmatrix} \begin{pmatrix}
0 \\ 0 \\ \beta^{-1} \\ 1
\end{pmatrix} \Bigg|_{S_0} = -\left(\dfrac{\alpha}{\alpha + \beta z} + \alpha  + \beta z\right).
\end{align*}}
A graph embedding $\phi(x,y,z,\varepsilon)$ of the {slow} manifold $S_0^{\varepsilon}$ is

\begin{align*}
\phi(x,y,z,\varepsilon) = \phi_0(x,y,z) + \ldots = \begin{pmatrix}
x \\ y \\ z \\ \frac{\alpha}{\alpha + \beta z} 
\end{pmatrix} + \cdots.
\end{align*}
Further, we have the following:

\begin{align*}
\Pi^{S_0} =  \begin{pmatrix}
1 & 0 & 0 & 0\\
0 & 1 & 0 & 0\\
0 & 0 &  \sigma_1 & -\beta^{-1} \sigma_1\\
0 & 0 & -\beta \sigma_2 & \sigma_2\\
\end{pmatrix},\quad
D\phi_0^L = \begin{pmatrix}
1 & 0 & 0 & 0\\
0 & 1 & 0 & 0\\
0 & 0 & 1 & 0
\end{pmatrix},
\end{align*}
where 
\begin{align*}
\sigma_1 = \dfrac{(\alpha + \beta z)^2}{(\alpha + \beta z)^2 + \alpha}, \quad \sigma_2 = \dfrac{\alpha}{(\alpha + \beta z)^2 + \alpha}.
\end{align*}
The reduced vector field \eqref{reduced:gen-param} is then given by:
\begin{align}
\begin{pmatrix}
\frac{dx}{dt} \\
\frac{dy}{dt} \\
\frac{dz}{dt} \\
\end{pmatrix} = \varepsilon \begin{pmatrix}
\frac{\alpha}{\alpha + \beta z} - \tilde{\rho}_2 x \\
\tilde{\rho}_3 x - \tilde{\rho}_4 y \\
 K(z) \left((\alpha+\beta z)(\tilde{\rho}_5 y - \tilde{\rho}_6 z)  - \tilde{\gamma} z\right)
\end{pmatrix}, \label{goodwin_eq3}
\end{align}
where 
\begin{align*}
K(z) = \dfrac{(\alpha + \beta z)}{(\alpha + \beta z)^2 + \alpha}.
\end{align*}

\begin{remark}
    {{This model has been previously derived for the case $\gamma = \rho_6$ and is called the \emph{pre-factor} QSSA in \cite{kimjosicbennett2,kimjosicbennett}}.}
\end{remark}

{Figure \ref{KF_discrepancy}} shows good agreement between the 4D model \eqref{KF-4D-nondim} {(solid blue lines)} and the 3D model reduction \eqref{goodwin_eq3} {(dashed green lines)}. Furthermore, recall from Remark \ref{KF_assumption} that the Kim-Forger model requires $\gamma = \rho_6$. The model is known for its oscillatory behavior; see e.g. \cite{kimtyson,kim2021} and references therein. However, {Figure \ref{KF_discrepancy}} shows that for the given set of parameter values, there is an onset of oscillations away from the Kim-Forger case $\gamma = \rho_6$ and this is correctly predicted by the model reduction \eqref{goodwin_eq3}. {The high accuracy of the reduction is evidenced by the dashed green trajectories overlaying the solid blue trajectories of the full system.}

\begin{figure}[ht]
\centering
\begin{subfigure}{0.49\textwidth}
\centering
 \includegraphics[width=1\linewidth]{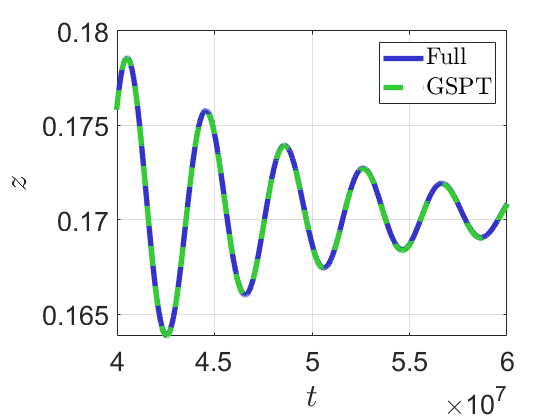}
 \subcaption{}
\end{subfigure}
\begin{subfigure}{0.49\textwidth}
\centering
 \includegraphics[width=1\linewidth]{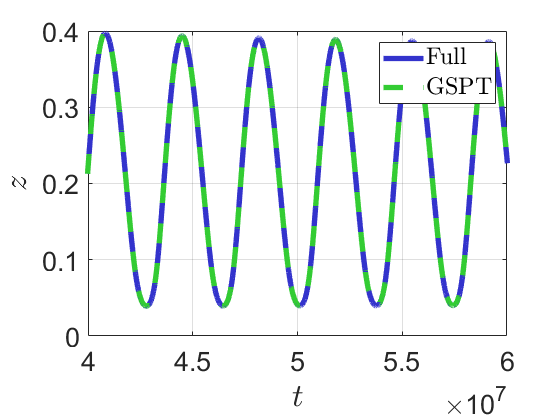}
 \subcaption{}
\end{subfigure}
\caption{{Numerical validation of the model reduction for the Kim-Forger (KF) reaction model.} A comparison of numerical simulations of the 4D full system \eqref{KF-4D-nondim} {(solid blue lines)} and the ci-GSPT model reduction \eqref{goodwin_eq3} {(dashed green lines)}. (a) $\gamma = \rho_6 = 10^{-6}$ {shows non-oscillatory dynamics}. (b) $\gamma = 1.5 \rho_6$, where $\rho_6 = 10^{-6}$, {shows oscillatory dynamics}. Other parameter values and ICs are given in Table \ref{KF_params}. {The dashed green reduction solution is visible on top of the solid blue full solution, demonstrating the high accuracy of the reduction.}}
 \label{KF_discrepancy}
\end{figure}

\begin{table}[ht]
    \centering
    \begin{tabular}{|c|c|c|c|}
    \hline
         Parameter/IC &  Value & Parameter/IC &  Value \\ \hline
         $x(0),y(0)$ & 1 & $\rho_2$ & $10^{-6}$ \\
         ${z}(0),s(0)$ & 0.06128 & $\rho_3$ & $ 
         10^{-5}$ \\
          $\alpha$ & $0.004$ & $\rho_4$ & $ 10^{-6}$ \\
         $\beta$ & $1$ & $\rho_5$ & $10^{-6}$ \\
          $\rho_1$ & $5\times 10^{-6}$ & $\rho_6$ & $10^{-6}$ \\
         \hline
    \end{tabular}
    \caption{Parameter values and ICs for simulations in Figure \ref{KF_discrepancy}. ICs are calculated so that the flow of the full 4D system \eqref{KF-4D-nondim} begins on the critical manifold. It is calculated in a similar way as in Remark \ref{remark_IC_tQSSA}.}
    \label{KF_params}
\end{table}

\section{Geometric Classification of the MM Reaction Scheme}
\label{MM_class}
In this section, we analyze the 67 singularly perturbed cases of the reversible MM reaction scheme (and the 23 cases for the irreversible scheme) to classify their geometric configurations. This analysis yields two classifications: one based on the geometry of the critical manifold, and another on the structure of the fast fibers, $N_0$. The latter provides a direct link between {the geometric framework outlined in Section \ref{sec:ciGSPTtools}} and the traditional QSSA approach.

 \subsubsection*{Critical manifold geometry}
{Consider the classical irreversible sQSSA case (where $\varepsilon:= \beta$, $\alpha,\gamma = \mathcal{O}(1), \delta = 0$). For this parameter configuration, system \eqref{eq:dimless-EK-2Dirr} becomes
\begin{align*}
    \begin{pmatrix}
        \frac{ds}{dt} \\ \frac{dc}{dt}
    \end{pmatrix} &= \begin{pmatrix}
0 \\ 1
\end{pmatrix} (s(1-c) - \alpha c - \gamma c)  + \varepsilon \Bigg( \begin{pmatrix}
        -s(1-c) \\ 0
    \end{pmatrix} + \alpha \begin{pmatrix}
        c \\ 0
    \end{pmatrix}\Bigg) \\
    &= N_0 f_0(s,c) + \varepsilon F_1(s,c).
\end{align*}
The set $S_0 = \{(s,c) \, | \, f_0(s,c) = 0 \}$ defines a one-dimensional manifold, which can be expressed as a graph over $s$:
\begin{align*}
c_0(s) = \dfrac{s}{\alpha + \gamma + s}.
\end{align*}
This is the well-known form of the critical manifold in the literature (see e.g. \cite{segelslemrod}). One may then calculate the non-trivial eigenvalue $Df_0 N_0 |_{S_0} = -(s+\alpha+\gamma)$; see Lemma \ref{non-trivial_eig}. Since it is negative, the critical manifold is attracting.} 

{
If we further assume that one of $\alpha$ or $\gamma$ is $\mathcal{O}(\varepsilon)$, the denominator of $c_0(s)$ retains the same functional form (dominated by the remaining $\mathcal{O}(1)$ parameter). We group these cases under the geometry `Form 1'.\\
However, if we assume that both $\alpha,\gamma = \mathcal{O}(\varepsilon)$, the manifold degenerates into two branches given by $S_{0,1} \cup S_{0,2} = \{s=0\}\cup\{c = 1\}$. We classify this geometry as `Form 2a'. It can then be shown that (i) $S_{0,1}$ is degenerate, i.e., the non-trivial eigenvalue is zero along $S_{0,1}$, (ii) $S_{0,2}$ is attracting for $ s>0$, and loses normal hyperbolicity at $s=0$, i.e., at $S_{0,1} \cap S_{0,2} = (0,1)$.}

{Now, consider a reversible sQSSA case (where $\varepsilon:= \beta$, $\alpha,\gamma,\delta = \mathcal{O}(1)$). For this parameter configuration, we have
\begin{align*}
    \begin{pmatrix}
        \frac{ds}{dt} \\ \frac{dc}{dt}
    \end{pmatrix} &=  \begin{pmatrix}
0 \\ 1
\end{pmatrix} (s(1-c) - (\alpha+\gamma) c + \delta (1-s)(1-c))
+ \varepsilon  \begin{pmatrix}
        -s(1-c) + \alpha c\\ -\delta c(1-c)
    \end{pmatrix} \\
    &= N_0 f_0(s,c) + \varepsilon F_1(s,c).
\end{align*}
The set $S_0 = \{(s,c) \, | \, f_0(s,c) = 0 \}$ is a one-dimensional manifold expressible as a graph over $s$:
\begin{align} \label{crit_man_beta_small}
c_0(s) = \dfrac{\delta + s - \delta s}{\alpha + \gamma + \delta + s- \delta s},
\end{align}}{which we designate as `Form 4'. We note that this function describes two disjoint critical manifold branches (hyperbolas) with horizontal asymptote $c = 1$ and vertical asymptote  $s=(\alpha+\gamma+\delta)/(\delta-1)$, $\delta\neq 1$,
and has also been derived in \cite{noethenwalcher,goekewalcherzerz}. Its upper branch is repelling and its lower branch is attracting; see Supplementary Material I for further details on its normal hyperbolicity.}

{
If we further assume that $\alpha,\gamma= \mathcal{O}(\varepsilon)$, the geometry of the critical manifold degenerates to $S_{0,1} \cup S_{0,2} =\left\{s=\frac{\delta}{\delta - 1}\right\} \cup \{c=1\}$. We classify this case also as `Form 2a', as it represents a shift of the vertical critical manifold branch observed in the irreversible case $\delta = 0$.}

{
To illustrate the interaction with asymptotically large parameters, consider a variation of the reversible sQSSA case (where $\varepsilon:= \beta$, $\alpha = \mathcal{O}(1)$, and $\gamma,\delta = \mathcal{O}\left(\frac{1}{\varepsilon}\right)$). In this regime, we obtain
\begin{align} 
c_0(s) = \dfrac{\tilde{\delta} - \tilde{\delta} s}{\tilde{\delta} + \tilde{\gamma} - \tilde{\delta} s}.
\end{align}
Similar to \eqref{crit_man_beta_small}, this defines hyperbolas with vertical and horizontal asymptotes. Hence, we also classify this case as `Form 4'.}

{For brevity, we do not explicitly show the derivations for all critical manifold geometries. We emphasize that the geometry is strictly determined by the zero level set of the leading-order term of system \eqref{eq:dimless-EK}; see Definition \ref{criticalmanifold_def}. The resulting geometric classification for all 67 singularly perturbed reversible cases is summarised in the following definition.}

\begin{definition}[Classification based on the geometry of critical manifolds] \label{criticalmanifold_def}
    For the 67 singularly perturbed cases, the following geometries arise; see Figure \ref{CMs_rev}:
    \begin{itemize}
        \item Form 1: $\left\{ c = \frac{s}{s+\Delta}\right\}$ where $\Delta \in \{\alpha,\gamma, \alpha+\gamma\}$,
        \item Form 2a: $\{s=s^*\}\cup\{c=1\}$ where $s^* \in \left\{0, 1 ,\frac{\delta}{\delta - 1} \right\}$
        \item Form 2b: $\{c=0\}$,
        \item Form 3: $\left\{ c = \frac{h_1(s) +\sqrt{h_2(s)}}{2 \beta \delta}\right\}\cup\left\{c = \frac{h_1(s) -\sqrt{h_2(s,c)}}{2 \beta \delta}\right\}$ where $h_1(s) = \delta + \gamma +\delta (\beta-s)$ and $h_2(s) = (\delta(\beta+s)-(\gamma + \delta))^2+4\beta \gamma \delta $.
        \item Form 4: $\left\{ c = \frac{\delta + g(s)}{\Delta + \delta + g(s)}\right\}$ where $g(s) = \{ -\delta s, s - \delta s\}$ and $\Delta \in \{\alpha,\gamma, \alpha+\gamma\}$,
        \item Form 5a: $\{c=-\beta^{-1}s + \beta^{-1}\}\cup\{c=1\}$,
        \item Form 5b: $\{ c = 0\}\cup\{ c = 1\}$,
        \item Form 5c: $\{ c = 1\}$.
    \end{itemize}
\end{definition}
\begin{remark}
    With a slight abuse of notation, we have dropped the tilde $\tilde{}$, recognizing that parameters in this definition are of order $\mathcal{O}(1)$; see Assumption \ref{parameterassumption}.
\end{remark}

\begin{figure}[ht]
\centering
 \includegraphics[width=0.58\linewidth]{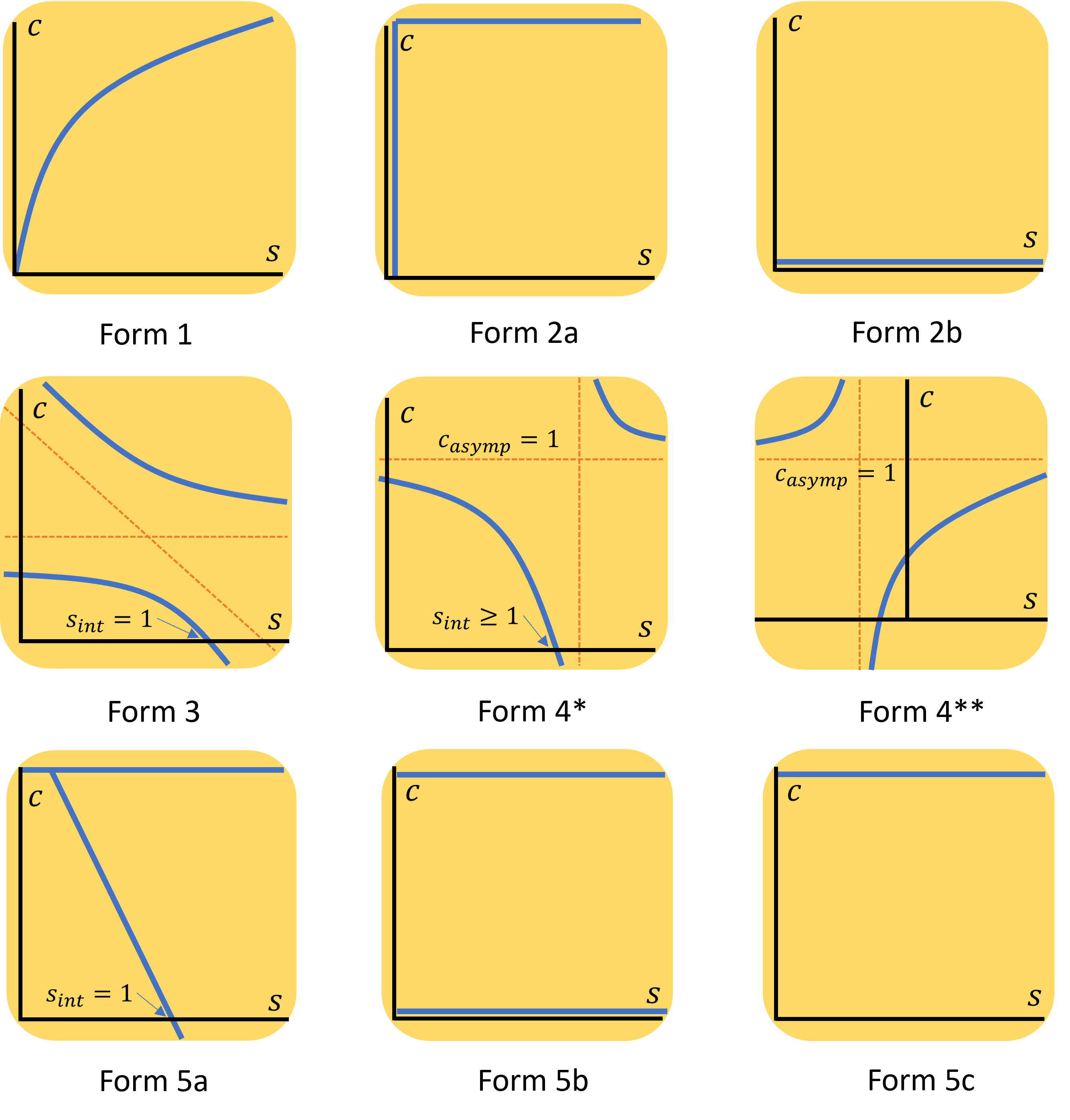}
\caption{Sketches of critical manifold branches (blue) arising from the MM reaction scheme, with some important features highlighted such as horizontal asymptotes and $s$-intercepts. For cases with Form 4 and $g(s) = s - \delta s$ (see Definition \ref{criticalmanifold_def}), $\delta > 1$ results in Form 4$^*$, $\delta < 1$ results in Form 4$^{**}$, and $\delta = 1$ is discussed in Supplementary Material I. Remaining cases where $g(s) = - \delta s$ result in Form 4$^*$ and $s_{int} = 1$.} \label{CMs_rev}
\end{figure}
{We note that in the 9 possible geometries, we either have (a) a unique
critical manifold branch (Forms 1, 2b, 5c), (b) two branches that intersect (Forms 2a, 5a), or (c) two branches that are disjoint (Forms 3, 4, 5b). This is sketched in Figure \ref{CMs_rev}.
Furthermore, only Forms 1, 2a (with vertical branch $\{s = s^\ast=0\}$) and 2b arise in the irreversible MM reaction $\delta=0$ and in fact, the irreversible cases are a subset of the reversible cases where $\delta = \mathcal{O}(\varepsilon)$.} 
{
\begin{remark}
    We make the following observations regarding the geometries listed in Definition \ref{criticalmanifold_def}:
    \begin{itemize}
        \item Forms 2a and 2b may be viewed as limiting cases of Form 1.
        \item Forms 5a, 5b, and 5c may be viewed as limiting cases of Form 3.
        \item Forms 3 and 4 represent hyperbolas, where Form 3 possesses an oblique asymptote.
    \end{itemize}
\end{remark}}
\bigskip
\subsubsection*{Fast fiber geometry}
\begin{definition}[Fast (linear) fiber classification] \label{fast_fiber_def}
  In system \eqref{eq:dimless-EK-2D}, we classify singular perturbation problems into Class S, R, or T, based on the orientation of $N_0$ as follows:
\begin{itemize}
    \item $N_0$ is vertical (Class S),
    \item $N_0$ is horizontal (Class R),
    \item $N_0$ is neither vertical nor horizontal (Class T). 
\end{itemize}
\end{definition}
{The relationship between the existing literature and our fast fiber classification is as follows: Class \textbf{S} corresponds to a parameter regime where the \textbf{S}tandard QSSA is applicable, Class \textbf{R} corresponds to where the \textbf{R}everse QSSA is applicable, and Class \textbf{T} corresponds to where the \textbf{T}otal QSSA is applicable. Recall that in the literature: (i) the sQSSA assumes that $c$ is fast (implying vertical fast fibers), (ii) the rQSSA assumes $s$ is fast (implying horizontal fast fibers), and (iii) the tQSSA utilizes {$\bar{s} = s+\beta c$} as the slow variable (implying transversal fast fibers); see, e.g., \cite{hta,segelslemrod,schnellmaini2000,schnellmaini2002,borghans,tzafriri}.} 

 {This allows us to define a {\em subclassification} based on the fast fiber orientation $N_0$ and the geometry of the critical manifold branches. We designate a case to be in Subclass S.1 if it possesses vertical fast fibers (S) and a critical manifold of Form 1. For example, consider the two cases described in Section \ref{sec:MM_GSPT}: (i) `the classical tQSSA' ($\varepsilon := \gamma$, $ \alpha,\beta = \mathcal{O}(1)$ and $\delta = 0$) belongs to Subclass T.1, and (ii) `an sQSSA case with multiple timescales' ($\varepsilon := \beta$, $\alpha = \mathcal{O}(1), \gamma = \mathcal{O}(\varepsilon)$ and $\delta = 0$) belongs to Subclass S.1. {Figure \ref{partition_rev} shows a schematic diagram combining fast fiber orientation and critical manifold geometry.}}

\begin{figure}[H]
\centering
\begin{subfigure}{1\textwidth}
\centering
 \includegraphics[width=0.83\linewidth]{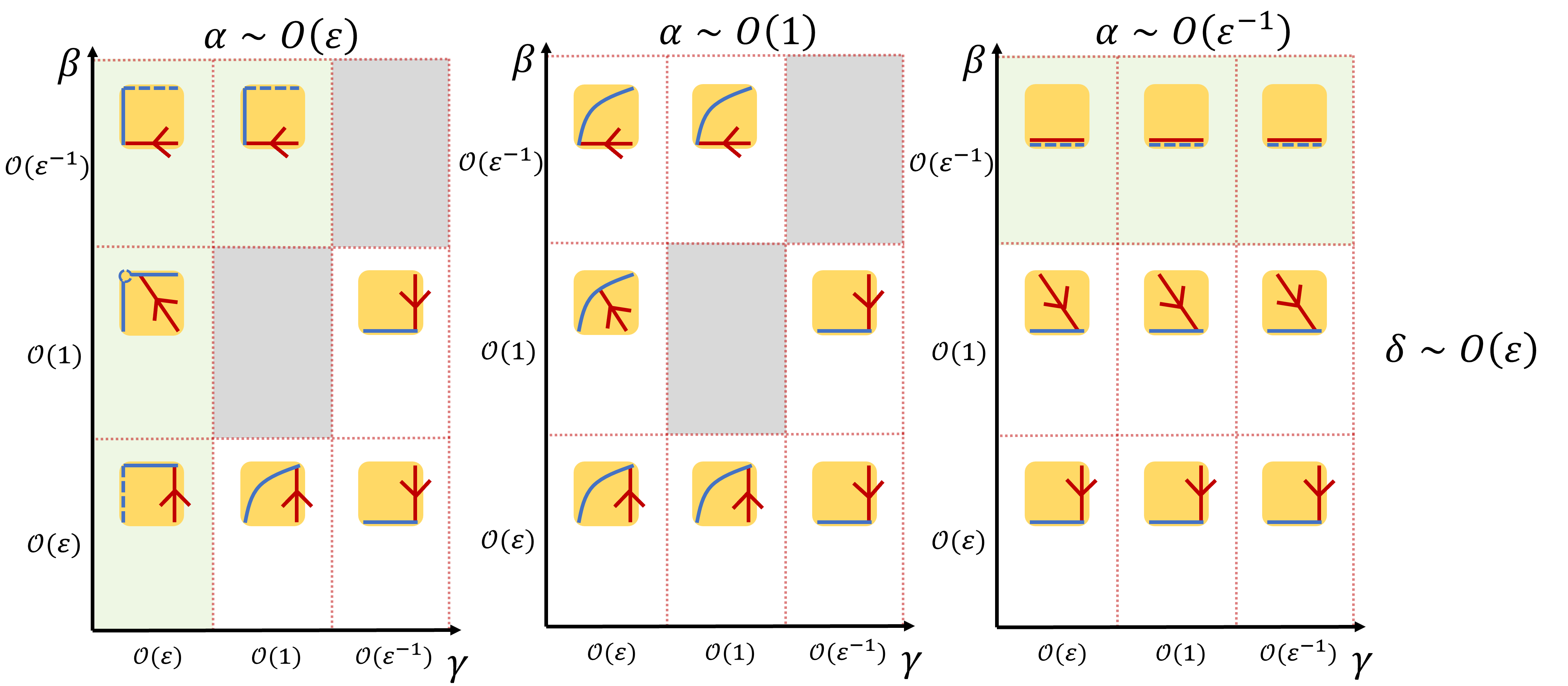}
 \centering
 \end{subfigure} 
\begin{subfigure}{1\textwidth}
\centering
  \includegraphics[width=0.83\linewidth]{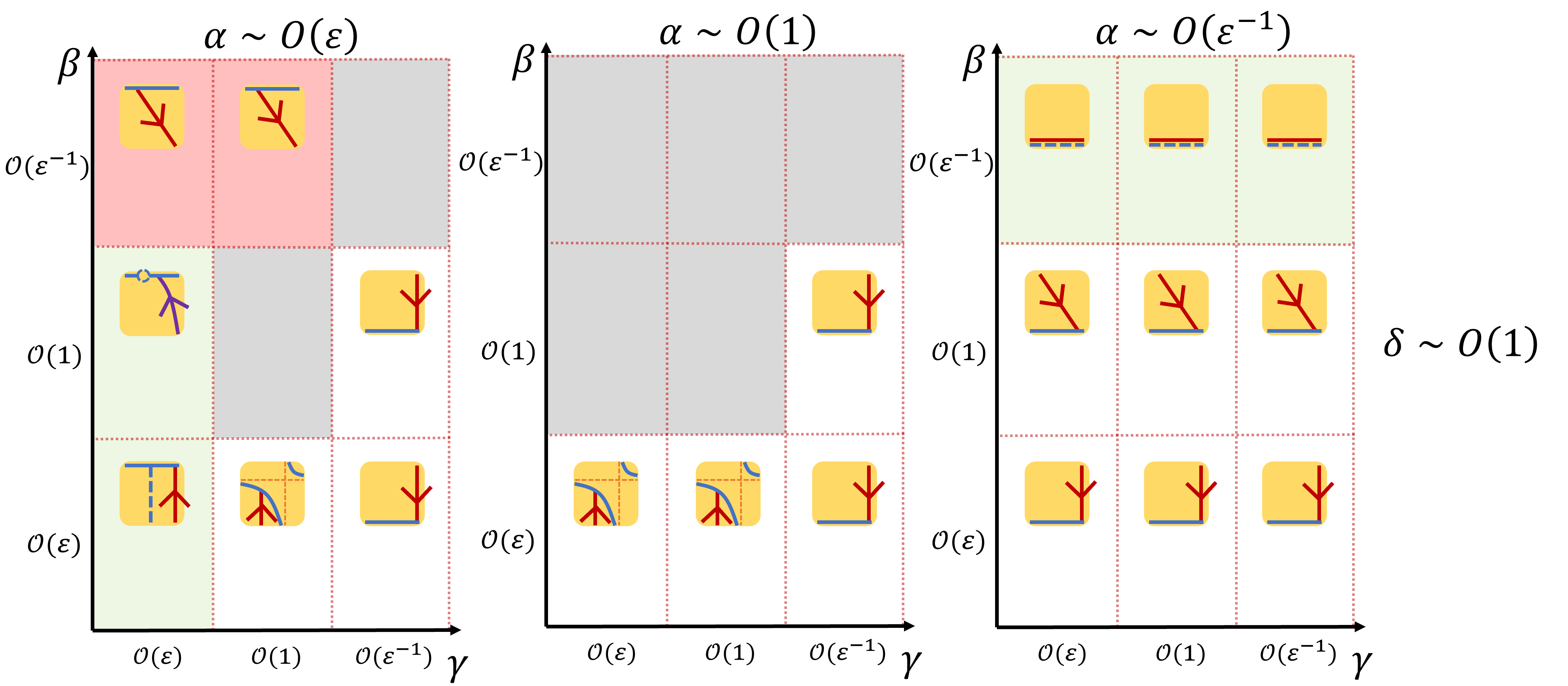}
\end{subfigure} 
\begin{subfigure}{1\textwidth}
\centering
  \includegraphics[width=0.83\linewidth]{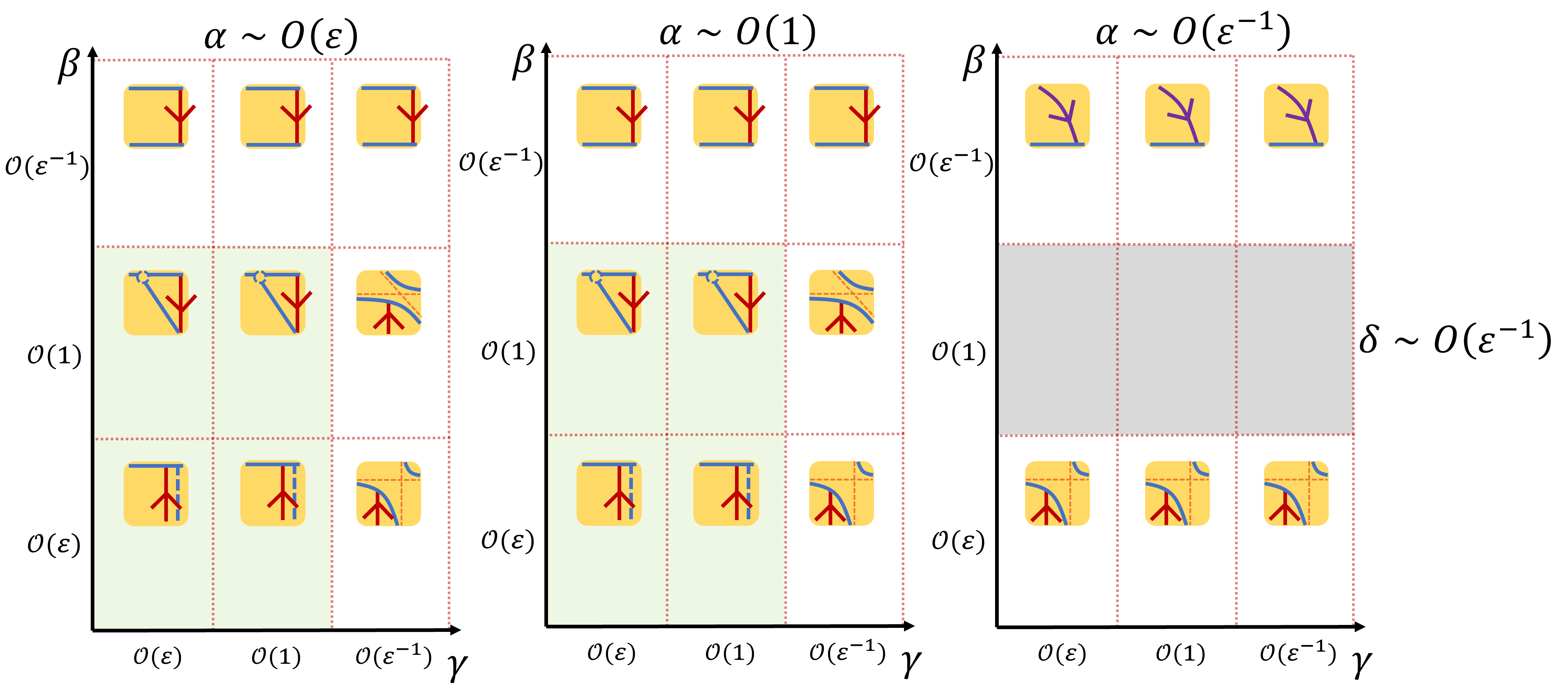}
\end{subfigure} 
\caption{A diagram of all 67 singularly perturbed MM cases, where $\delta = \mathcal{O}(\varepsilon)$ shows also the irreversible cases. Critical manifold(s) $S_0$ (blue) and a fast fiber $\mathcal{W}^S(y)$ where $y \in S_0$ (red: linear, purple: non-linear) are shown. 
Gray boxes: non-singularly perturbed cases. {Green boxes}: the critical manifold has degeneracy. Dotted lines and circled points denote location. Red boxes: cases with a single critical manifold and it is repelling. Supplementary Material I provides more information on parameter configuration. 
}  \label{partition_rev}
\end{figure}

\begin{remark} \label{remark_straightfibers}
    In system \eqref{eq:dimless-EK-2D}, 63 out of 67 singularly perturbed reversible cases have $\mathcal{W}^S = \mathcal{N}_0$; see also Remark \ref{remark_IC_tQSSA}. The 4 exceptions are: one case in Subclass T.5c where $\alpha,\gamma = \mathcal{O}(\varepsilon)$,$\beta,\delta = \mathcal{O}(1)$, and three cases in Subclass T.2b where $\alpha,\beta,\gamma = \mathcal{O}(\varepsilon^{-1})$ and $\delta$ is either $\mathcal{O}(\varepsilon), \mathcal{O}(1)$ or $\mathcal{O}(\varepsilon^{-1})$; see Figure \ref{partition_rev}. In such cases, the classical coordinate transformation $\bar{s} = s + \beta c$ does not globally result in vertical fast fibers.
\end{remark}
{
\begin{remark} \label{no_need}
    Irrespective of the subclass, the {tools} outlined in Section \ref{sec:ciGSPTtools} provides a  {framework for calculating a specific model reduction for a given parameter configuration.}
\end{remark}}

\section{Irreversible MM: Literature Comparison and Model Reductions}
\label{lit_comp_irr}
Tables \ref{sQSSA_class}, \ref{tQSSA_class}, and \ref{rQSSA_class} classify the 23 singularly perturbed irreversible cases into their respective subclasses. {Table \ref{comparison_epsilon} defines various small parameters ($\varepsilon_{HTA},\varepsilon_{SS}, \varepsilon_{SS,SM}, 
\varepsilon_{BdBS}, \varepsilon_{T}$) derived in the literature and compares these 23 cases to selected QSSA validity conditions.} This effectively aligns our fast fiber classification (see Definition \ref{fast_fiber_def}) with established validity criteria.

We emphasize that the classical approach in the literature is to broaden the validity regions of QSSAs, either for a specific type such as the sQSSA\footnote{Classically, $\varepsilon_{HTA} \ll 1$ is the condition for the sQSSA, but \cite{segelslemrod} present the condition $\varepsilon_{SS} \ll 1$, which encapsulates more parameter configurations.} or for QSSAs in general\footnote{For example, \cite{borghans,tzafriri} state that their tQSSA validity conditions overlap with and extend those of the sQSSA and rQSSA. This expands the parameter region where a QSSA is valid.}. Consequently, it is expected that a more recent condition for the sQSSA, such as $\varepsilon_{SS} \ll 1$, also covers cases that fall into our Class T. Indeed, the two main conditions for tQSSA validity, $\varepsilon_{BdBS} \ll 1$ and $\varepsilon_T \ll 1$, individually encapsulate 22 of the 23 singularly perturbed cases we have enumerated.
{

\begin{table}[H]
\centering
\begin{tabular}{ |c|c|c| } 
 \hline
Subclass S.1  & Subclass S.2a & Subclass S.2b \\ \hline
i. $\beta = \mathcal{O}(\varepsilon) , \alpha, \gamma = \mathcal{O}(1)$ & i. $\beta ,\alpha, \gamma = \mathcal{O}(\varepsilon)$ & i. $\beta  = \mathcal{O}(\varepsilon), \alpha = \mathcal{O}\left(\frac{1}{\varepsilon}\right), \gamma = \mathcal{O}(1)$ \\
ii. $\beta,\alpha = \mathcal{O}(\varepsilon), \gamma = \mathcal{O}(1)$ & & ii. $\beta  = \mathcal{O}(\varepsilon), \alpha = \mathcal{O}(1), \gamma = \mathcal{O}\left(\frac{1}{\varepsilon}\right)$\\
iii. $\beta,\gamma = \mathcal{O}(\varepsilon), \alpha = \mathcal{O}(1)$ & & iii. $\beta  = \mathcal{O}(\varepsilon), \alpha,\gamma = \mathcal{O}\left(\frac{1}{\varepsilon}\right)$  \\ 
& & iv. $\beta,\gamma  = \mathcal{O}(\varepsilon), \alpha = \mathcal{O}\left(\frac{1}{\varepsilon}\right)$ \\
& & v. $\beta,\alpha = \mathcal{O}(\varepsilon), \gamma = \mathcal{O}\left(\frac{1}{\varepsilon}\right)$ \\ 
& & vi. $\beta,\alpha = \mathcal{O}(1)$, $\gamma = \mathcal{O}\left(\frac{1}{\varepsilon}\right)$ \\
& & vii. $\beta = \mathcal{O}(1), \alpha = \mathcal{O}(\varepsilon), \gamma = \mathcal{O}\left(\frac{1}{\varepsilon}\right)$ \\
\hline
\end{tabular} 
\caption{The 11 cases in Class S for the irreversible MM reaction scheme.}\label{sQSSA_class}
\end{table}
\begin{table}[ht]
\centering
\begin{tabular}{ |c|c|c| }
 \hline
Subclass T.1 & Subclass T.2a & Subclass T.2b \\  \hline
i. $\gamma = \mathcal{O}(\varepsilon), \alpha, \beta = \mathcal{O}(1) $ & i. $\gamma,\alpha = \mathcal{O}(\varepsilon), \beta = \mathcal{O}(1)$ & i. $\alpha = \mathcal{O}\left(\frac{1}{\varepsilon}\right), \beta, \gamma = \mathcal{O}(1)$ \\ 
& & ii. $\alpha = \mathcal{O}\left(\frac{1}{\varepsilon}\right), \beta = \mathcal{O}(1)$, \\
& & $\gamma = \mathcal{O}(\varepsilon)$ \\
& & iii. $\alpha,\gamma = \mathcal{O}\left(\frac{1}{\varepsilon}\right), \beta = \mathcal{O}(1)$
\\ \hline
\end{tabular} 
\caption{The 5 cases in Class T for the irreversible MM reaction scheme.}\label{tQSSA_class}
\end{table}
\begin{table}[ht]
\centering
\begin{tabular}{ |c|c| }
 \hline
Subclass R.1 & Subclass R.2a 
 \\ \hline
 i. $\beta = \mathcal{O}\left(\frac{1}{\varepsilon}\right), \alpha, \gamma = \mathcal{O}(1)$ & i. $\beta = \mathcal{O}\left(\frac{1}{\varepsilon}\right), \alpha = \mathcal{O}(\varepsilon), \gamma = \mathcal{O}(1)$\\
 ii. $\beta = \mathcal{O}\left(\frac{1}{\varepsilon}\right), \alpha = \mathcal{O}(1), \gamma = \mathcal{O}(\varepsilon)$  &  ii. $\beta = \mathcal{O}\left(\frac{1}{\varepsilon}\right), \alpha, \gamma = \mathcal{O}(\varepsilon)$ \\ \hline
\multicolumn{2}{|c|}{Subclass R.2b } \\ \hline
\multicolumn{2}{|c|}{i.  $\beta,\alpha = \mathcal{O}\left(\frac{1}{\varepsilon}\right), \gamma = \mathcal{O}(1)$ } \\
\multicolumn{2}{|c|}{ii. $\beta,\alpha = \mathcal{O}\left(\frac{1}{\varepsilon}\right), \gamma = \mathcal{O}(\varepsilon)$} \\
\multicolumn{2}{|c|}{iii. $\beta,\alpha,\gamma = \mathcal{O}\left(\frac{1}{\varepsilon}\right)$}
\\ \hline
\end{tabular}
\caption{The 7 cases in Class R for the irreversible MM reaction scheme.}\label{rQSSA_class}
\end{table}

\begin{table}[H]
    \centering
    \begin{tabular}{|m{4.9cm}|p{0.8cm}|c|} \hline
       Condition and QSSA type & Class & Cases  \\ \hline
       $\varepsilon_{HTA}:= \beta \ll 1$  & S & 1.i,ii,iii, 2a.i, 2b.i,ii,iii,iv,v \\ 
       sQSSA & R & \\
       \cite{hta} & T & \\ \hline
       $\varepsilon_{SS}:=\frac{\beta}{1+\alpha+\gamma} \ll 1$   & S & 1.i,ii,iii, 2a.i, 2b.i,ii,iii,iv,v,vi,vii \\ 
       sQSSA & R & \\
      \cite{segelslemrod} & T & 2b.i,ii,iii \\ \hline
       $\varepsilon_{HTA}^{-1} \ll 1, \varepsilon_{SS,SM} : = \frac{\gamma}{\beta} \ll 1$ & S & \\ 
       rQSSA & R & 1.i,ii, 2a.i,ii, 2b.i,ii  \\ 
       \cite{schnellmaini2000, segelslemrod} & T & \\ \hline
       $\varepsilon_{SS}^{-1} \ll 1$  & S &  \\ 
       rQSSA & R & 1.i,ii, 2a.i,ii \\
       \cite{schnellmaini2000} & T & \\ \hline
       $\varepsilon_{BdBS} := \frac{\beta \gamma}{(\alpha+\gamma+\beta+1)^2} \ll 1$ & S & 1.i,ii,iii, 2a.i, 2b.i,ii,iii,iv,v,vi,vii   \\ 
       tQSSA & R & 1.i,ii, 2a.i,ii, 2b.i,ii \\
       \cite{borghans} & T & 1.i, 2a.i, 2b.i,ii, 2b.iii \\ \hline
       $\varepsilon_{T}:= \frac{\gamma}{2}((1-\frac{4}{\gamma} \varepsilon_{BdBS})^{\frac{-1}{2}}-1) \ll 1$  & S & 1.i,ii,iii, 2a.i, 2b.i,ii,iii,iv,v,vi,vii   \\ 
       tQSSA & R & 1.i,ii, 2a.i,ii, 2b.i,ii \\
       \cite{tzafriri} & T & 1.i, 2a.i, 2b.i,ii,iii \\ 
       \hline
    \end{tabular}
    \caption{The distribution of the 23 irreversible cases based on selected validity conditions for sQSSA, rQSSA and tQSSA in literature. The parameter configurations for each case are shown in Tables \ref{sQSSA_class}, \ref{tQSSA_class} and \ref{rQSSA_class}. }
    \label{comparison_epsilon}
\end{table}

\begin{remark}
    The Reich-Sel'kov condition $\frac{\beta}{\alpha + \gamma}$ \cite{reichselkov} is also frequently employed to assess sQSSA validity. For small $e_0$, this expression is proportional to the leading-order eigenvalue ratio of the equilibrium of system \eqref{eq:dimless-EK-2Dirr} at the origin (see, e.g., \cite{accuracy} and references therein). We emphasize that a significant spectral gap (timescale separation) is a prerequisite for the ci-GSPT tools presented here (see Definition \ref{def_normallyhyperbolic}). Similarly, the condition $\varepsilon_{CSP} = \frac{\beta \gamma (1-c) }{(\alpha+\gamma + s + \beta (1-c))^2}$, which also approximates an eigenvalue ratio, is utilized in \cite{mmvalid} to estimate timescale disparity globally across the $(s,c)$-phase space.
\end{remark}}

In contrast, {the} ci-GSPT approach provides a constructive method for model reduction, i.e., given a specific parameter configuration, the steps in Section \ref{sec:ciGSPTtools} can be followed to derive model reductions. This process yields a unique model reduction for each normally hyperbolic critical manifold of any singularly perturbed case.

We also note that our fast fiber classification leads to an interesting observation: the assumption on $\beta$ alone {\em almost exclusively} dictates which class a singularly perturbed case falls in, with the exception of 2 cases in Class S where $\beta = \mathcal{O}(1)$; see Figure \ref{venn_irreversible}.

Table \ref{summary_table_ST} presents the model reductions for 14 (out of 16) S/T cases, each possessing a unique, normally hyperbolic, and attracting critical manifold. Similarly, Table \ref{summary_table_ST3} provides the model reductions for 2 additional cases, derived for their critical manifolds away from any degeneracies. The remaining singularly perturbed cases are addressed in the concluding remarks of the section.

\begin{remark}
    Cases exhibiting a loss of normal hyperbolicity can be resolved via the blow-up method  \cite{dumortierroussarie1996,kuehn2015,wechselberger2020}, which magnifies the dynamics in the vicinity of non-hyperbolic points, revealing local behavior that is inaccessible via the ci-GSPT reduction in Section \ref{sec:ciGSPTtools}. Lax and Walcher \cite{laxwalcher} provide general conditions under which the rescaled system retains a singularly perturbed structure, thereby allowing for the continued application of the presented ci-GSPT tools. For Case S.2a.i, {the rescaling $s \to \varepsilon \tilde{s}$ is required, which is termed by Lax and Walcher as a `degenerate scaling' and is analogous to the procedure described in Remark \ref{rescale_time} involving degenerate scalings of the system parameters.}
\end{remark}}

\begin{remark}[Class R] \label{rQSSA_discard}
    All cases in Subclass R.1 and R.2a undergo a rapid equilibriation. This means that the reduction is trivial, as the layer flow terminates on a stable equilibrium of the model reduction, which is $s^* = 0$ for cases in Subclass R.1 and $c^*= 0$ for cases in Subclass R.2a. The equilibrium persists as a stable node at the origin of the 2D system \eqref{eq:dimless-EK-2Dirr}; see Remark \ref{eq_remark}  and {Supplementary Material III}. We note that the model reductions of these cases are not shown in this paper but may be calculated using the tools presented. Furthermore, since cases in Subclass R.2b have a single critical manifold that is degenerate everywhere, a blow-up (rescaling) is required for further analysis.
\end{remark}

\begin{remark} \label{eq_irr}
    All cases in Table \ref{summary_table_ST}  have a stable equilibrium at $s^* = 0$ and Subclass T.2a in Table \ref{summary_table_ST3} has a stable equilibrium at $c^* = 0$. This persists as a stable node at the origin of the 2D system \eqref{eq:dimless-EK-2Dirr}; see Remark \ref{eq_remark}  and {Supplementary Material III}. 
\end{remark}

\begin{remark}
    For Subclass T.2a, there are two reductions: (a) on $S_{0,1} = \{s=0\}$ and (b)  $S_{0,2} = \{c=1\}$. However, note $S_{0,1}$ cannot be written as a graph over $s$ and so we have done the reduction over the $c$-coordinate. Figure 3 in \cite{algorithmic2023} is an example where `switching' between sQSSA and rQSSA is required. This is encapsulated in the two model reductions for Case T.2a.i shown in Table \ref{summary_table_ST3}. 
\end{remark}

\begin{remark}
    Since we {apply a consistent framework} to each individual case -- as opposed to relying on the overarching reductions often found in the literature -- one may observe `apparent' discrepancies between the ci-GSPT reduction and the {model reductions and product formation rates derived via standard QSSAs.} However, Goeke et. al. \cite{goekewalcherzerz} demonstrates that QSSAs and singular perturbation reductions typically agree to leading-order in the small parameter. Furthermore, Supplementary Material II provides a detailed comparison between the product formation rates derived from our ci-GSPT approach and those in the literature.
\end{remark}

{
\begin{remark}
 The ci-GSPT tools presented in this paper aim to provide a solution to the misuse of QSSAs {outside of their validity regimes} highlighted by Kim and Tyson \cite{kimtyson}. We address this by outlining {an approach} that constructs specific model reductions appropriate for a given parameter regime. However, Kim and Tyson also highlight potential misuse in the stochastic setting. For instance, the validity of the sQSSA in a stochastic context requires robustness to initial conditions away from the slow manifold (induced by stochastic fluctuations) \cite{kimjosicbennett,kimtyson}. From a geometric perspective, this suggests that the correct projection of initial conditions onto the slow manifold (along the fast fibers) is crucial; see e.g. \cite{eilertsenstroberg}.
\end{remark}}

\begin{figure}[ht]
\centering
\includegraphics[scale=0.45]{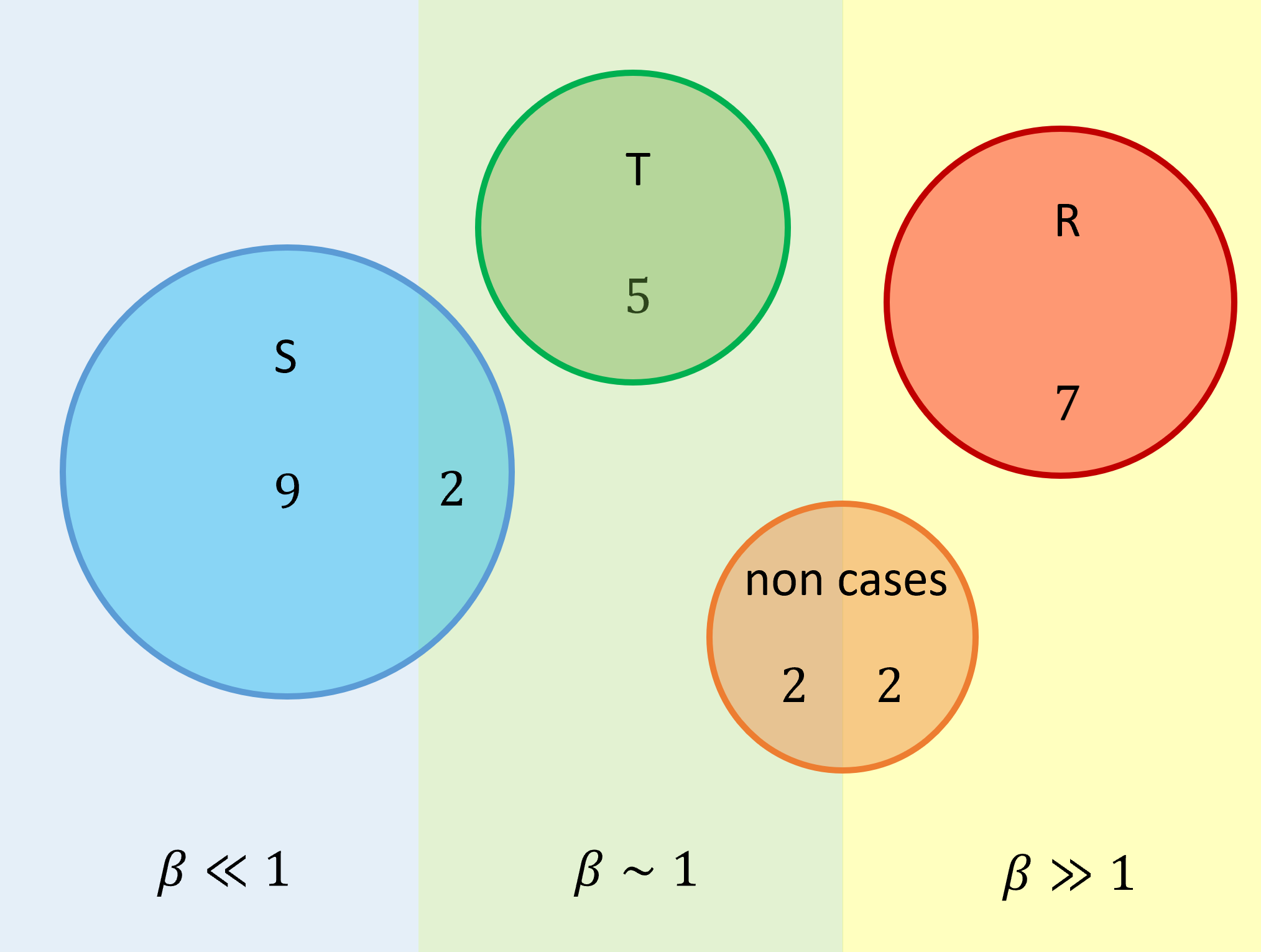}
    \caption{A Venn diagram of the distribution of cases for each class and its assumption on $\beta$.}
\label{venn_irreversible}
\end{figure}

\begin{table}[H]
    \centering
    \begin{tabular}{|m{1.4cm}|m{1.6cm}|m{1.85cm}|m{1.4cm}|m{4.5cm}|}
    \hline
        Subclass (cases) & Parameter assumptions $(\alpha,\beta,\gamma)$ order & Substrate depletion & Product formation  & Further comments \\ \hline
         S.1 (i,ii) & $(-,\varepsilon,1)$ & $- \varepsilon \frac{\gamma s}{\Delta + s}$ & $\varepsilon \frac{\gamma s}{\Delta  + s}$ & $\varepsilon := \beta$. (i) $\alpha = \mathcal{O}(1)$, $\Delta = \alpha+\gamma$. (ii) $\alpha = \mathcal{O}(\varepsilon)$, $\Delta = \gamma$.  \\[2.5ex]
         S.1 (iii) & $(1,\varepsilon,\varepsilon)$&  $-\varepsilon^2 \frac{\tilde{\gamma} s}{\alpha + s}$  & $\varepsilon^2 \frac{\tilde{\gamma} s}{\alpha + s}$ & $\varepsilon := \beta$.  \\[2.5ex]
        S.2b (i,iv) & $\left(\frac{1}{\varepsilon},\varepsilon,-\right)$ & $-\varepsilon^n \frac{\gamma}{\tilde{\alpha}} s$ & $\varepsilon^n \frac{\gamma}{\tilde{\alpha}} s$ & (i) $\gamma = \mathcal{O}(1)$, $n=3$. (ii) {$\gamma = \mathcal{O}(\varepsilon) \implies n=4$}. \\[2.5ex]
        S.2b (ii,iii,v) & $(-,\varepsilon,1)$ & $-\varepsilon^2 \Delta s$ & $\varepsilon^2  \Delta s $ & $\varepsilon := \beta$.   (iii) $\alpha = \mathcal{O}\left(\frac{1}{\varepsilon}\right)$, $\Delta = \{ \frac{\tilde{\gamma}}{\tilde{\alpha} + \tilde{\gamma}}\}$. (ii) $\alpha = \mathcal{O}(1)$, (v) $\alpha = \mathcal{O}(\varepsilon)$ with  $\Delta = 1$ for both. \\[2.5ex]
        S.2b (vi, vii) &  $\left(-,1,\frac{1}{\varepsilon}\right)$ & $-\varepsilon \beta s$ & $\varepsilon \beta s$  & $\varepsilon := \frac{1}{\gamma}$. (vi) $\alpha = \mathcal{O}(1)$. (vii) $\alpha = \mathcal{O}(\varepsilon)$. \\[2.5ex]
        T.1 (i) & $(1,1,\varepsilon)$ &  $-\varepsilon \frac{ \beta s ( \alpha  +s)}{\alpha \beta + (\alpha + s)^2} $ &  $\varepsilon \frac{\beta s}{s+\alpha}$ & $\varepsilon := \gamma$. \\[2.5ex]
        T.2b (i, ii) & $\left(\frac{1}{\varepsilon},1,-\right)$ & $-\varepsilon^n \beta \gamma s$ & $\varepsilon^n \beta \gamma s$ & $\varepsilon := \frac{1}{\alpha}$. (i) $\gamma = \mathcal{O}(1)$, $n = 2$. (ii){$\gamma = \mathcal{O}(\varepsilon) \implies n=3$}. \\[2.5ex]
        T.2b (iii) & $\left(\frac{1}{\varepsilon},1,\frac{1}{\varepsilon}\right)$ & $-\varepsilon \beta \tilde{\gamma} \frac{s}{1+\tilde{\gamma}}$ & $\varepsilon \beta \tilde{\gamma} \frac{s}{1+\tilde{\gamma}}$ & $\varepsilon := \frac{1}{\alpha}$. \\ [0.75ex]
        \hline
    \end{tabular}
    \caption{ Summary of the model reductions of 14 cases in Class S and T of the irreversible MM. The $\tilde{}$ notation indicates if the parameter has been scaled; see Assumption \ref{parameterassumption}. The $-$ notation under the second column indicates that the parameter assumption is under `Further comments'.}
    \label{summary_table_ST}
\end{table}

\begin{table}[ht]
    \centering
    \begin{tabular}
    {|m{1.3cm}|m{1.6cm}|m{1.8cm}|m{1.5cm}|m{4.6cm}|}
    \hline
        Subclass (cases) & Parameter assumptions $(\alpha,\beta,\gamma)$ order & Substrate depletion & Product formation & Further comments \\ \hline
          S.2a (i) & $(\varepsilon,\varepsilon,\varepsilon)$ & $-\varepsilon^2 \tilde{\gamma}$ & $\varepsilon^2 \tilde{\gamma}$ & $\varepsilon := \beta$. Results shown for $S_{0,2} = \{c=1\}$. Loss of normal hyperbolicity on $S_{0,1} = \{s=0\}.$ \\[2.5ex]
        T.2a (i) & $(\varepsilon,1,\varepsilon)$  & (a) $\dot{c} = -\varepsilon c$, (b) $\dot{s} = -\varepsilon \beta $  & (a) $\varepsilon \beta c$, (b) $\varepsilon \beta $ & $\varepsilon := \gamma$. (a) $S_{0,1} = \{s=0\}$; critical manifold not a graph over $s$. Product formation: $\dot{p} = - \beta \dot{c}$. (b) $S_{0,2} = \{c=1\}$. Loss of normal hyperbolicity at $(0,1)$. \\
        \hline
    \end{tabular}
    \caption{Summary of the model reductions of the 2 cases in Class S and T with a loss of normal hyperbolicity of the irreversible MM. The $\tilde{}$ notation indicates if the parameter has been scaled; see Assumption \ref{parameterassumption}.}
    \label{summary_table_ST3}
\end{table}

{

\section{Reversible MM: Relevant Cases and Model Reductions}
\label{sec:rev}

While the reversible MM reaction scheme ($\delta>0$) has received less attention in the literature, it has frequently served as a key example for the development of ci-GSPT tools (see, e.g., \cite{noethenwalcher,goekewalcher2013,goekewalcher}). In this section, we derive the corresponding model reductions. Given the large number of possibilities -- 67 singularly perturbed cases in total (listed in Supplementary Material I) -- we begin by setting aside cases that are either not biologically relevant or that lead to trivial reductions identifiable without significant calculation.

\subsubsection*{Negative leading-order product formation rate conditions}
Biologically, it is desirable to create more product $p$ over time, and so we discard those cases with negative leading-order product formation rates. 
The product formation rate\footnote{The product formation rate used here is the last equation of \eqref{eq:dimless-EK}.} is
\begin{align}
    \dfrac{dp}{dt} &= \beta \gamma c(s,\varepsilon) - \beta \delta (1-\beta c(s,\varepsilon) - s)(1-c(s,\varepsilon)) \nonumber \\
    &\approx \beta \gamma c_0(s) - \beta \delta (1-\beta c_0(s) - s)(1-c_0(s))\,,  \label{product_formation}
\end{align}
where $c(s,\varepsilon) \approx c_0(s)$ is the leading-order approximation $S_{0}$ of $S_{0}^\varepsilon$. This leads to a negative leading-order product formation rate if
\begin{align}
   (\beta = \mathcal{O}(\varepsilon) \,\,\, \text{and} \,\,\,   \gamma \ll \delta) \quad \quad \text{or}  \quad \quad (\beta = \mathcal{O}(1) \,\,\, \text{and} \,\,\,   \gamma \ll \delta). \label{condition1}
\end{align}
Furthermore, if $c_0(s) \equiv 0$, as in Form 2b and the bottom branch of Form 5b, we assume the leading-order approximation $c(s,\varepsilon) \approx \varepsilon c_1(s)$\footnote{We cannot know without further calculations if $c_1 \equiv 0$ also, due to the iterative nature of the parametrization method; see Section \ref{sec:ciGSPTtools}.} and obtain a negative leading-order product formation rate if:
{\begin{align}
   &(\beta = \mathcal{O}(\varepsilon) \,\,\, \text{and} \,\,\, \gamma = \mathcal{O}(\delta)) \quad \quad \text{or}  \quad \quad (\beta = \mathcal{O}(1) \,\,\,\text{and} \,\,\, \gamma = \mathcal{O}(\delta))  \nonumber \\
   &\text{or} \quad \quad \left(\beta = \mathcal{O}(\varepsilon^{-1}) \,\,\, \text{and} \,\,\,  (\gamma = \mathcal{O}(\delta) \,\,\, \text{or} \,\,\, \gamma \ll \delta)\right). \label{condition2}
\end{align}}
\subsubsection*{Biologically relevant and non-trivial cases}
 Tables \ref{sQSSA_class_rev_pfr} and \ref{tQSSA_class_rev_pfr} lists the 27 remaining biologically relevant (i.e., does not satisfy \eqref{condition1} or \eqref{condition2}, and at least (a subset of) one critical manifold is attracting; see Remark \ref{discard}) and non-trivial cases (i.e. the model reduction is non-trivial; see Remark \ref{rQSSA_discard_rev}).

\begin{remark} \label{discard}
The following comments are on discarded cases based on biological relevance. The reader should also refer to Figure \ref{partition_rev} for the schematic diagram of the subclassification and Supplementary Material I for the corresponding parameter configurations.

\begin{itemize}
    \item We discard cases if 
    \eqref{condition1} or \eqref{condition2} are satisfied. In particular:
    \begin{itemize}
\item Since cases in Subclass S.4b have two critical manifolds $S_{0,1} = \{c=0\}$ and $S_{0,2} = \{c=1\}$, these have been discarded based on the lower branch $S_{0,1}$. 
    \item  Cases S.2a.iiii,iv,v and vi are discarded based on its normally hyperbolic horizontal branch $S_{0,2} = \{c=1\}$, as the result is a negative product formation rate. However, we should also note that the IC $(1,0)$ is on the degenerate vertical branch $S_{0,1} = \{s=1\}$.
\end{itemize}
\item Cases T.5c.ii and iii are discarded as its single critical manifold is repelling.
\end{itemize}

\end{remark}

\begin{remark}[Class R; c.f. Remark \ref{rQSSA_discard}] \label{rQSSA_discard_rev}
        All reversible cases in Subclass R.1 and R.2a have the same parameter configurations as in the irreversible case in addition to $\delta = \mathcal{O}(\varepsilon)$; c.f. Table \ref{rQSSA_class} and Table 3 in Supplementary Material I. An additional $\mathcal{O}(\varepsilon)$-term then slightly perturbs the stable node at the origin to $0 < s^*,c^* < 1$ ({Supplementary Material III} proves the bounds of the equilibrium value). Hence, such reversible cases undergo a {rapid equilibriation}. Indeed, one can check that the {model reductions of cases in Subclass R.1 have a stable equilibrium at $s^* = 0$ and the model reductions of cases in Subclass R.2a have a stable equilibrium at $c^*=0$}. Furthermore, since Subclass R.2b has a unique branch that is degenerate everywhere, we cannot discard it based on the sufficient conditions in \eqref{condition1} and \eqref{condition2}. However, as we did in the irreversible case, we do not consider it further in this paper as it requires a blow-up (rescaling).
\end{remark}

\begin{table}[ht]
\centering
\begin{tabular}{ |c|c| } 
 \hline
Subclass S.1 & Subclass S.2a  
 \\ \hline 
 i. $\beta, \delta = \mathcal{O}(\varepsilon), \alpha, \gamma = \mathcal{O}(1)$ & i. $\beta, \alpha, \gamma, \delta = \mathcal{O}(\varepsilon)$
 \\
 ii. $\beta, \alpha, \delta = \mathcal{O}(\varepsilon), \gamma = \mathcal{O}(1)$ & \\
 iii. $\beta,  \gamma, \delta = \mathcal{O}(\varepsilon),\alpha = \mathcal{O}(1)$ & \\ \hline
Subclass S.2b & Subclass S.3 \\ \hline
i. $\beta, \alpha, \delta = \mathcal{O}(\varepsilon), \gamma = \mathcal{O}\left(\frac{1}{\varepsilon} \right)$  & i. $\beta = \mathcal{O}(1), \alpha = \mathcal{O}(\varepsilon), \gamma, \delta = \mathcal{O}\left( \frac{1}{\varepsilon}\right)$ \\
ii. $\beta, \alpha, = \mathcal{O}(\varepsilon), \gamma = \mathcal{O}\left(\frac{1}{\varepsilon} \right),\delta = \mathcal{O}(1)$ &    ii. $\beta, \alpha = \mathcal{O}(1), \gamma, \delta = \mathcal{O}\left( \frac{1}{\varepsilon}\right)$ \\ \cline{2-2}
iii. $\beta,\delta = \mathcal{O}(\varepsilon), \alpha = \mathcal{O}(1), \gamma = \mathcal{O}\left(\frac{1}{\varepsilon}\right)$   &  Subclass S.4 
\\  \cline{2-2}
iv. $\beta = \mathcal{O}(\varepsilon), \alpha, \delta = \mathcal{O}(1), \gamma = \mathcal{O}\left(\frac{1}{\varepsilon}\right)$ & i.$\beta, \alpha, = \mathcal{O}(\varepsilon), \gamma,\delta = \mathcal{O}\left(\frac{1}{\varepsilon} \right)$  
\\
 vii.  $\beta, \delta = \mathcal{O}(\varepsilon), \alpha = \mathcal{O}\left(\frac{1}{\varepsilon}\right), \gamma = \mathcal{O}(1)$ &  ii. $\beta, \alpha, = \mathcal{O}(\varepsilon), \gamma, \delta = \mathcal{O}(1) $ 
\\
 ix. {$\beta,\delta = \mathcal{O}(\varepsilon), \alpha, \gamma = \mathcal{O}\left(\frac{1}{\varepsilon}\right)$}  &  iv. $\beta = \mathcal{O}(\varepsilon), \alpha, \gamma, \delta = \mathcal{O}(1)$ 
\\
x.  $\beta = \mathcal{O}(\varepsilon), \alpha, \gamma = \mathcal{O}\left(\frac{1}{\varepsilon}\right), \delta = \mathcal{O}(1)$ & v. $\beta = \mathcal{O}(\varepsilon), \alpha = \mathcal{O}(1), \delta, \gamma = \mathcal{O}\left(\frac{1}{\varepsilon}\right)$ 
\\
xi. $\beta= \mathcal{O}(1), \alpha, \delta = \mathcal{O}(\varepsilon), \gamma = \mathcal{O}\left( \frac{1}{\varepsilon}\right)$  & viii. $\beta = \mathcal{O}(\varepsilon), \alpha, \gamma, \delta = \mathcal{O}\left(\frac{1}{\varepsilon}\right)$ 
\\ 
xii. $\beta, \delta = \mathcal{O}(1), \alpha = \mathcal{O}(\varepsilon), \gamma = \mathcal{O}\left( \frac{1}{\varepsilon}\right)$ &
\\
xiii. $\beta, \alpha = \mathcal{O}(1), \gamma = \mathcal{O}\left( \frac{1}{\varepsilon}\right) , \delta = \mathcal{O}(\varepsilon)$ & 
\\
 xiv. $\beta, \alpha,\delta = \mathcal{O}(1), \gamma = \mathcal{O}\left( \frac{1}{\varepsilon}\right)$ & 
\\
\hline
\end{tabular}
\caption{The 22 relevant cases in Class S for the reversible MM reaction scheme. }\label{sQSSA_class_rev_pfr}
\end{table}

\begin{table}[ht]
\centering
\begin{tabular}{ |c|c| } 
 \hline
Subclass T.1 & Subclass T.2a  
 \\ \hline
i. $\gamma,\delta = \mathcal{O}(\varepsilon), \alpha, \beta = \mathcal{O}(1)$  & i. $\gamma,\alpha ,\delta = \mathcal{O}(\varepsilon), \beta = \mathcal{O}(1)$  \\ \hline
 Subclass T.2b &  \\ \hline
 iii. $\alpha = \mathcal{O}\left(\frac{1}{\varepsilon}\right), \beta,\gamma = \mathcal{O}(1), \delta = \mathcal{O}(\varepsilon)$ & \\ 
  v. $\alpha,\gamma = \mathcal{O}\left(\frac{1}{\varepsilon}\right), \beta = \mathcal{O}(1), \delta = \mathcal{O}(\varepsilon)$ & \\
  vi. $\alpha, \gamma = \mathcal{O}\left(\frac{1}{\varepsilon}\right), \beta, \delta = \mathcal{O}(1)$ & \\
\hline
\end{tabular}
\caption{The 5 relevant cases in Class T for the reversible MM reaction scheme.}\label{tQSSA_class_rev_pfr}
\end{table}

Tables \ref{summary_table_ST_rev}, \ref{summary_table_ST_rev_new}, \ref{summary_table_ST2_rev1}, \ref{summary_table_ST2_rev} and \ref{summary_table_ST3_rev} show a model reduction for 25 cases with at least one critical manifold that is normally hyperbolic and attracting. Table \ref{loss_rev} shows model reductions away from degeneracies for 2 cases. We note that Subclass S.1, S.2b, T.1, T.2b, S.2a and T.2a (see Tables \ref{summary_table_ST_rev}, \ref{summary_table_ST_rev_new}, \ref{summary_table_ST3_rev} and \ref{loss_rev}) are also subclasses in the irreversible MM (see Tables \ref{summary_table_ST} and \ref{summary_table_ST3}).

  \begin{remark}[c.f. Remark \ref{eq_irr}]
    All cases in Tables \ref{summary_table_ST_rev}, \ref{summary_table_ST_rev_new}, \ref{summary_table_ST2_rev1}, \ref{summary_table_ST2_rev} and \ref{summary_table_ST3_rev} have a stable equilibrium for $0 \leq s^* < 1$ and Case T.2a.i in Table \ref{loss_rev} has a stable equilibrium for $0 < c^* < 1$. This persists as a stable node $(s^*,c^*)$, where $0 < s^*,c^* < 1$, in the 2D system \eqref{eq:dimless-EK-2D}; see Remark \ref{eq_remark} and Supplementary Material III. The reader can find the proofs that the {relevant} equilibrium for Cases T.1.i, T.2a.i, S.3.i and S.3.ii is in the correct range and is attracting in Supplementary Material III, as the calculations are more involved.
\end{remark}

  \begin{table}[H]
    \centering
    \begin{tabular}{|m{1.6cm}|m{1.6cm}|m{2.3cm}|m{2.1cm}|m{3.2cm}|}
    \hline 
        Subclass (cases) & Parameter assumptions $(\alpha,\beta,\gamma,\delta)$ order & Substrate depletion & Product formation & Further comments \\ \hline
        S.1 (i, ii) & $(-,\varepsilon,1,\varepsilon)$ & $-\varepsilon \frac{\gamma s}{\Delta + s}$ & $\varepsilon \frac{\gamma s}{\Delta + s}$ &$\varepsilon := \beta$. (i) $\alpha = \mathcal{O}(1)$, $\Delta = \alpha+\gamma$. (ii) $\alpha = \mathcal{O}(\varepsilon)$, $\Delta = \gamma$.  \\[2.5ex]
         S.1 (iii) &  $(1,\varepsilon,\varepsilon,\varepsilon)$ & $-\varepsilon^2 \frac{\tilde{\gamma}s +  \alpha \tilde{\delta} (s-1)}{\alpha + s} $ & $\varepsilon^2 \frac{\tilde{\gamma}s +  \alpha \tilde{\delta} (s-1)}{\alpha + s}$ & $\varepsilon := \beta$. \\ [1ex]
         \hline
    \end{tabular}
    \caption{Summary of 3 relevant cases in Subclass S.1 with a single attracting critical manifold of the reversible MM. Subclass S.1 also appear in the irreversible MM; see Table \ref{summary_table_ST} for comparison. The $\tilde{}$ notation indicates if the parameter has been scaled; see Assumption \ref{parameterassumption}. The $-$ notation under the second column indicates that the parameter assumption is under `Further comments'.}
    \label{summary_table_ST_rev}
\end{table}

 \begin{table}[H]
    \centering
    \begin{tabular}{|m{1.6cm}|m{1.6cm}|m{2.3cm}|m{2.1cm}|m{3.2cm}|}
    \hline 
        Subclass (cases) & Parameter assumptions $(\alpha,\beta,\gamma,\delta)$ order & Substrate depletion & Product formation & Further comments \\ \hline
S.2b (i,ii,iii, iv) & $\left(-,\varepsilon,\frac{1}{\varepsilon},-\right)$ & $-\varepsilon^2 s$ &  $\varepsilon^2 s$ & $\varepsilon := \beta$. (i) $\alpha,\delta = \mathcal{O}(\varepsilon)$. (ii) $\alpha = \mathcal{O}(\varepsilon), \delta = \mathcal{O}(1)$. (iii) $\alpha = \mathcal{O}(1), \delta = \mathcal{O}(\varepsilon)$. (iv)  $\alpha,\gamma = \mathcal{O}(1)$. \\ [2.5ex]
        S.2b (xi,xii,xiii, xiv) & $\left(-,1,\frac{1}{\varepsilon},-\right)$ & $-\varepsilon \beta s$ &  $\varepsilon \beta s$ &  $\varepsilon := \frac{1}{\gamma}$. (xi) $\alpha, {\delta}  = \mathcal{O}(\varepsilon)$. 
        (xii)  $\alpha = \mathcal{O}(\varepsilon), {\delta}  = \mathcal{O}(1)$. (xiii) $\alpha = \mathcal{O}(1), {\delta}  = \mathcal{O}(\varepsilon)$. (xiv) $\alpha,{\delta}  = \mathcal{O}(1)$. \\ [2.5ex] 
        S.2b (ix) & $\left(\frac{1}{\varepsilon},\varepsilon,\frac{1}{\varepsilon},\varepsilon\right)$ & $-\varepsilon^2 \frac{\tilde{\gamma}s}{\tilde{\alpha} + \tilde{\gamma}}$ &  $-\varepsilon^2 \frac{\tilde{\gamma}s}{\tilde{\alpha + \tilde{\gamma}}}$  &  $\varepsilon := \beta$. Note $c_1(s) \not \equiv 0$.\\ [2.5ex] 
        S.2b (x) &  $\left(\frac{1}{\varepsilon},\varepsilon,\frac{1}{\varepsilon},1\right)$  & $-\varepsilon^2  \frac{\tilde{\alpha} (-\delta - s + \delta s)}{\tilde{\alpha} + \tilde{\gamma}}   -\varepsilon^2 s$ &  $\varepsilon^2  \frac{\tilde{\alpha} (-\delta - s + \delta s)}{\tilde{\alpha} + \tilde{\gamma}} + \varepsilon^2 s$  &  $\varepsilon := \beta$. Note $c_1(s) \not \equiv 0$. \\ [2.5ex]
        S.2b (vii) &   $\left(\frac{1}{\varepsilon},\varepsilon,1,\varepsilon\right)$  & $-\varepsilon^3 \tilde{\delta} (\tilde{s}-1)- \varepsilon^3 \frac{\gamma s}{\tilde{\alpha}}$  & $\varepsilon^3 \tilde{\delta} (\tilde{s}-1) + \varepsilon^3 \frac{\gamma s}{\tilde{\alpha}}$ & $\varepsilon := \beta$ and $\delta = \mathcal{O}(\varepsilon)$ $\alpha= \mathcal{O}\left( \frac{1}{\varepsilon}\right)$. Note $c_1(s) \not \equiv 0$.\\ 
         \hline
    \end{tabular}
    \caption{Summary of 11 relevant cases in S.2b with a single attracting critical manifold of the reversible MM. Subclass S.2b also appear in the irreversible MM; see Table \ref{summary_table_ST} for comparison. The $\tilde{}$ notation indicates if the parameter has been scaled; see Assumption \ref{parameterassumption}. The $-$ notation under the second column indicates that the parameter assumption is under `Further comments'.}
    \label{summary_table_ST_rev_new}
\end{table}

\begin{table}[ht]
    \centering
    \begin{tabular}{|m{1.3cm}|m{1.6cm}|m{2.4cm}|m{2.2cm}|m{3.3cm}|}
    \hline 
        Subclass (cases) & Parameter assumptions $(\alpha,\beta,\gamma,\delta)$ order & Substrate depletion & Product formation & Further comments \\ \hline
         S.3 (i, ii) &  $\left(-,1,\frac{1}{\varepsilon},\frac{1}{\varepsilon}\right)$  & $-\varepsilon \frac{H(s)}{2 \tilde{\delta}} \,\,\, $ & $\varepsilon (\beta c'_0(s) +1) \frac{H(s)}{2 \tilde{\delta}}$ & $\varepsilon := \frac{1}{\gamma}$. (i) $\alpha = \mathcal{O}(\varepsilon)$. (ii) $\alpha = \mathcal{O}(1)$. \\ [1ex] 
        \hline
    \end{tabular}
     \caption{Summary of 2 relevant cases in Subclass S.3 with two critical manifolds. The reductions shown are for the attracting critical manifold. {$H(s) = h_3(s) - w(s) (1 -\sqrt{h_2(s)})$. (i) $h_3(s) = \tilde{\delta} s (\beta + s - 1), w(s) = s$ and (ii) $h_3(s) = - \tilde{\delta} (\alpha + s + \alpha \beta - \alpha s - \beta s - s^2)$, $w(s) = \alpha + s$. For their critical manifold and $h_2(s)$, see Definition \ref{criticalmanifold_def}, where the attracting critical manifold has the negative square root.}
     The $\tilde{}$ notation indicates if the parameter has been scaled; see Assumption \ref{parameterassumption}.  The $-$ notation under the second column indicates that the parameter assumption is under `Further comments'.}
    \label{summary_table_ST2_rev1}
\end{table}
\begin{table}[ht]
    \centering
    \begin{tabular}{|m{1.3cm}|m{1.6cm}|m{2.4cm}|m{2.2cm}|m{3.3cm}|}
    \hline 
        Subclass (cases) & Parameter assumptions $(\alpha,\beta,\gamma,\delta)$ order & Substrate depletion & Product formation & Further comments \\ \hline
        S.4 (i, ii) &  $(\varepsilon,\varepsilon,-,-)$  & $-\varepsilon^n\frac{\gamma s}{\gamma + \delta + g(s)}$ &  $\varepsilon^n\frac{\gamma s}{\gamma + \delta + g(s)}$ & $\varepsilon := \beta$. and $\alpha = \mathcal{O}(\varepsilon)$. (i): $\gamma,\delta = \mathcal{O}\left( \frac{1}{\varepsilon} \right)$, $n=2$ and $g(s) = -\tilde{\delta} s$. (ii) $\gamma,\delta = \mathcal{O}(1)$, $n=1$ and $g(s) = s-\delta s $. \\ [2.5ex]
        S.4 (iv, v) &  $(1,\varepsilon,-,-)$  & $-\varepsilon^n \frac{ (\gamma s +\alpha \delta (s-1))}{\gamma + \delta + h(s)}$ & $\varepsilon^n \frac{ (\gamma s +\alpha \delta (s-1))}{\gamma + \delta + h(s)}$ & $\varepsilon: = \beta$. (iv) $\gamma,\delta = \mathcal{O}(1)$, $n=1$ and $h(s) = \alpha + s - \delta s$. (v) $\gamma,\delta = \mathcal{O}\left( \frac{1}{\varepsilon} \right)$, $n=2$ and $h(s) = -\tilde{\delta} s$.  \\ [2.5ex]
        S.4 (viii) &  $\left(\frac{1}{\varepsilon},\varepsilon,\frac{1}{\varepsilon},\frac{1}{\varepsilon}\right)$ & $-\varepsilon \frac{ \tilde{\alpha} \tilde{\delta} (s-1)}{\tilde{\alpha} + \tilde{\gamma} + \tilde{\delta} - \tilde{\delta} s}$ & $\varepsilon \frac{ \tilde{\alpha} \tilde{\delta} (s-1)}{\tilde{\alpha} + \tilde{\gamma} + \tilde{\delta} - \tilde{\delta} s}$ & $\varepsilon: = \beta$. Rapid equilibriation; {Remark \ref{rQSSA_discard}}. \\ 
        \hline
    \end{tabular}
     \caption{Summary of 5 relevant cases in Subclass S.4 with two critical manifolds. The reductions shown are for the attracting critical manifold. The $\tilde{}$ notation indicates if the parameter has been scaled; see Assumption \ref{parameterassumption}.  The $-$ notation under the second column indicates that the parameter assumption is under `Further comments'.}
    \label{summary_table_ST2_rev}
\end{table}

 \begin{table}[H]
    \centering
    \begin{tabular}{|m{1.3cm}|m{1.6cm}|m{2.5cm}|m{2.2cm}|m{3.2cm}|}
    \hline 
        Subclass (cases) & Parameter assumptions $(\alpha,\beta,\gamma,\delta)$ order & Substrate depletion & Product formation & Further comments \\ \hline
        T.1 (i) &  $(1,1,\varepsilon,\varepsilon)$ &  $-\varepsilon H(s)$ & $\varepsilon (\beta c'_0(s) +1) H(s)$  &  $\varepsilon:= \gamma$. Critical manifold $c_0(s) = \frac{s}{\alpha + s}$. \\ [2.5ex] 
        T.2b (iii) & $\left(\frac{1}{\varepsilon},1,1,\varepsilon\right)$ &  $-\varepsilon^2 \beta ( \tilde{\delta}(s-1) + \gamma s)$ & $\varepsilon^2 \beta ( \tilde{\delta}(s-1) + \gamma s)$ & $\varepsilon:= \frac{1}{\alpha}$. $c_1(s) \not \equiv 0$. \\ [2.5ex]
        T.2b (v) & $\left(\frac{1}{\varepsilon},1,\frac{1}{\varepsilon},\varepsilon\right)$ &  $-\varepsilon \frac{\beta \tilde{\gamma s}}{1+\tilde{\gamma}}$ & $\varepsilon \frac{\beta \tilde{\gamma s}}{1+\tilde{\gamma}}$ & $\varepsilon:= \frac{1}{\alpha}$. $c_1(s) \not \equiv 0$. \\ [2.5ex]
        T.2b (vi) & $\left(\frac{1}{\varepsilon},1,\frac{1}{\varepsilon},1\right)$ &  $-\varepsilon \frac{\beta (-\delta + (\delta-1)s)}{1+\tilde{\gamma}} - \varepsilon \beta s$ &  $\varepsilon \frac{\beta (-\delta + (\delta-1)s)}{1+\tilde{\gamma}} + \varepsilon \beta s$  & $\varepsilon:= \frac{1}{\alpha}$. $c_1(s) \not \equiv 0$. \\
        \hline
    \end{tabular}
    \caption{Summary of 4 relevant cases in Class T with a single attracting critical manifold of the reversible MM. Subclass T.1 and T.2b also appear in the irreversible MM; see Table \ref{summary_table_ST} for comparison. {$H(s) = \frac{\beta \left( d_2 s^2 + d_1 s + d_0 \right)}{(\alpha+s)^2 + \alpha \beta}$ where $d_0 = - \alpha^2 \tilde{\delta}, d_1 = \alpha (1-\tilde{\delta} + \alpha \tilde{\delta} + \beta \tilde{\delta}), d_2 = 1+\alpha \tilde{\delta}$.} The $\tilde{}$ notation indicates if the parameter has been scaled; see Assumption \ref{parameterassumption}.}
    \label{summary_table_ST3_rev}
\end{table}

\begin{table}[H]
    \centering
    \begin{tabular}{|m{1.3cm}|m{1.6cm}|m{2cm}|m{1.9cm}|m{4.0cm}|}
    \hline
        Subclas (cases) & Parameter assumptions $(\alpha,\beta,\gamma,\delta)$& Substrate depletion & Product formation & Further comments \\ \hline
          S.2a (i) &  $(\varepsilon,\varepsilon,\varepsilon,\varepsilon)$  & $-\varepsilon^2 \tilde{\gamma}$ & $\varepsilon^2 \tilde{\gamma}$ & $\varepsilon := \beta$. Results for $S_{0,2} = \{c=1\}$. Loss of normal hyperbolicity on $S_{0,1} = \{s=0\}.$  \\[2.5ex]
        T.2a (i) & $(\varepsilon,1,\varepsilon,\varepsilon)$ & (a) {$\dot{c} = -\varepsilon (c -\tilde{\delta}(1-\beta c)(1-c))$}, (b) $\dot{s} = -\varepsilon \beta $  & (a) $\varepsilon \beta (c - \tilde{\delta} (1-\beta c)(1-c))$, (b) $\varepsilon \beta $ & $\varepsilon := \gamma$. Product formation: $\dot{p} = - \beta \dot{c}$.  (a) $S_{0,1} = \{s=0\}$; critical manifold not a graph over $s$. (b) $S_{0,2} = \{c=1\}$. Loss of normal hyperbolicity at $(0,1)$. \\ \hline
    \end{tabular}
    \caption{Summary of the 2 cases a loss of normal hyperbolicity in Class S and T of the reversible MM. Subclass S.2a and T.2a also appear in the irreversible MM; see Table \ref{summary_table_ST3} for comparison. The $\tilde{}$ notation indicates if the parameter has been scaled; see Assumption \ref{parameterassumption}. } \label{loss_rev}
\end{table}

\section{Conclusions}
\label{sec:conclusion}
Model reductions for chemical reaction networks (CRNs) are highly sought after, but common methods like the quasi-steady-state approximation (QSSA) present challenges {such as the ambiguity of selecting between simultaneously valid QSSAs for a given parameter regime.} In this paper, we approached model reduction differently. Using a combination of coordinate-independent geometric singular perturbation theory (ci-GSPT) and the parametrization method, {we calculated appropriate model reductions for a given parameter regime. This approach circumvents the need to select and validate a particular QSSA \textit{a priori}.}

Specifically, we demonstrated {the systematic application of these tools by} producing a comprehensive catalog of 14 and 25 distinct model reductions for the irreversible and reversible Michaelis-Menten reaction schemes, respectively (Tables \ref{summary_table_ST}, \ref{summary_table_ST_rev}, \ref{summary_table_ST_rev_new}, \ref{summary_table_ST2_rev1}, \ref{summary_table_ST2_rev} and \ref{summary_table_ST3_rev}). We also detailed cases leading to trivial reductions or a loss of normal hyperbolicity (Remarks \ref{rQSSA_discard}, \ref{rQSSA_discard_rev}; Tables \ref{summary_table_ST3}, \ref{loss_rev}) {and the dynamics of the degenerate subsets are the subject of ongoing investigation}. Furthermore, we demonstrated the scalability of {the} approach by {calculating a model reduction} for the more complex Kim-Forger oscillator.

In advocating for the proper use of QSSAs, Kim and Tyson \cite{kimtyson} note difficulties such as the complexity of the resulting model reduction and the challenge of finding a correct coordinate transformation. While {the use of ci-GSPT tools may also yield algebraically complex reductions, a key advantage of this approach} is that no coordinate transformations are required, which simplifies the comparison of dynamics between the full system and its reduction.

The parametrization method was essential for complementing ci-GSPT, allowing us to systematically compute higher-order approximations of the slow manifold and the slow dynamics defined on it. However, this method relies on the normal hyperbolicity of the critical manifold. Cases where this property is lost require the blow-up method for a complete resolution. {As noted in \cite{multiple}, symbolic calculations of higher-order derivatives may become computationally intensive, particularly for high-dimensional systems, and suggests the need for future work on algorithmic efficiency, such as the implementation of automatic differentiation \cite{automatic}. We note that as the parametrization method is well-established, such improvements have already been explored in various other contexts {to address computational complexity}; see e.g. \cite{param, roberts,robertsetal,huguetetal}.}

In summary, we demonstrated {the systematic use of ci-GSPT in conjunction with the parametrization method for calculating model reduction of CRNs}. While further work is required to fully analyze degenerate cases where loss of normal hyperbolicity occurs, our application to a large number of parameter configurations for benchmark systems provides a robust and versatile {approach} for analyzing complex biological systems.

\section*{Acknowledgments}
The authors thank Dimitris Goussis (Khalifa University, UAE), Vivien Kirk and James Sneyd (University of Auckland, NZ) for fruitful discussions about this manuscript. TEFL thanks, in particular, Vivien Kirk for travel support to the University of Auckland where part of this research has been conducted. Finally, we also thank Jae-Kyoung Kim and Eu-Min Jeong (KAIST, Korea) for discussions on the Kim-Forger model, Anthony Roberts (University of Adelaide, Australia) for further context regarding the parametrization method and the referees for their valuable comments. ChatGPT was used to proofread the text.

\bibliographystyle{siam}
\bibliography{references}

@Misc{siam,
  key = {zzz},
  title =	 {{SIAM} Style Manual: For journals and books},
  year =	 2013,
url = {https://epubs.siam.org/pb-assets/files/SIAM_STYLE_GUIDE_2019-1635349464967.pdf}}

@Misc{amsmath,
  author =	 {{American Mathematical Society}},
  title =	 {User's Guide for the \texttt{amsmath} Package
                  (Version 2.0)},
  url =		 {ftp://ftp.ams.org/pub/tex/doc/amsmath/amsldoc.pdf},
  urldate =	 {2015-07-30},
  year =	 2002}

@article{michaelismenten,
   author = {Michaelis, L. and Menten, M. L.},
   title = {{Die Kinetik der Invertinwirkung}},
   journal = {Biochem. Z.},
   volume = {49},
   pages = {333-369},
   year = {1913}
}

@book{henri,
   author = {Henri, Victor},
   title = {Lois générales de l'action des diastases},
   publisher = {Librairie Scientifique A. Hermann},
   year = {1903},
}

@article{briggshaldane,
   author = {Briggs, G. E. and Haldane, J. B. S.},
   title = {A note on the kinetics of enzyme action},
   journal = {Biochem. J.},
   volume = {19},
   number = {2},
   pages = {338-339},
   year = {1925},
}

@article{hta,
   author = {Heineken, F. G. and Tsuchiya, H. M. and Aris, R.},
   title = {On the mathematical status of the pseudo-steady state hypothesis of biochemical kinetics},
   journal = {Math. Biosci.},
   volume = {1},
   number = {1},
   pages = {95-113},
   ISSN = {0025-5564},
   year = {1967},
}

@article{segelslemrod,
   author = {Segel, L. A. and Slemrod, M.},
   title = {The quasi-steady-state assumption: a case study in perturbation},
   journal = {SIAM Rev.},
   volume = {31},
   number = {3},
   pages = {446-477},
   ISSN = {0036-1445},
   year = {1989},
}

@article{tzafriri,
   author = {Tzafriri, A. R.},
   title = {{Michaelis-Menten Kinetics at High Enzyme Concentrations}},
   journal = {Bull. Math. Biol.},
   volume = {65},
   number = {6},
   pages = {1111-1129},
   ISSN = {1522-9602},
   year = {2003},
}

@article{borghans,
   author = {Borghans, J. A. M. and De Boer, R. J. and Segel, L. A.},
   title = {Extending the quasi-steady state approximation by changing variables},
   journal = {Bull. Math. Biol.},
   volume = {58},
   number = {1},
   pages = {43-63},
   ISSN = {0092-8240},
   year = {1996},
}

@article{tikhonov,
   author = {Tikhonov, Andrei Nikolaevich},
   title = {Systems of differential equations containing small parameters in the derivatives},
   journal = {Mat. Sb.},
   volume = {73},
   number = {3},
   pages = {575-586},
   ISSN = {0368-8666},
   year = {1952},
}

@article{fenichel,
title = {Geometric singular perturbation theory for ordinary differential equations},
journal = {J. Differ. Equ.},
volume = {31},
number = {1},
pages = {53-98},
year = {1979},
issn = {0022-0396},
author = {N. Fenichel}
}

@article{goekewalcher,
   author = {Goeke, Alexandra and Walcher, Sebastian},
   title = {A constructive approach to quasi-steady state reductions},
   journal = {J. Math. Chem.},
   volume = {52},
   number = {10},
   pages = {2596-2626},
   ISSN = {1572-8897},
   year = {2014},
}

@article{goekewalcherzerz,
   author = {Goeke, Alexandra and Walcher, Sebastian and Zerz, Eva},
   title = {Classical quasi-steady state reduction—A mathematical characterization},
   journal = {Phys. D},
   volume = {345},
   pages = {11-26},
   ISSN = {0167-2789},
   year = {2017},
}

@InProceedings{goekewalcher2013,
  title={{Quasi-Steady State: Searching for and Utilizing Small Parameters}},
  author={Goeke, Alexandra
and Walcher, Sebastian},
  booktitle={{Recent
 Trends in Dynamical Systems.}},
editor={Johann, Andreas
and Kruse, Hans-Peter
and Rupp, Florian
and Schmitz, Stephan},
  pages={153-178},
  year={2013},
  series={Springer Proceedings of Mathematics \& Statistics},
volume = {35},
chapter = {8},
publisher = {Springer},
address = {Basel}}

@book{wechselberger2020,
  title={Geometric singular perturbation theory beyond the standard form},
  author={Wechselberger, M.},
  year={2020},
  publisher={Springer},
series = {Frontiers in Applied Dynamical Systems: Reviews and Tutorials},
volume = {6}
}

@article{accuracy,
   author = {Srivastava, Kashvi and Eilertsen, Justin and Booth, Victoria and Schnell, Santiago},
   title = {{Accuracy Versus Predominance: Reassessing the Validity of the Quasi-Steady-State Approximation}},
   journal = {Bull. Math. Biol.},
   volume = {87},
   number = {6},
   pages = {73},
   ISSN = {1522-9602},
   year = {2025},
}

@inbook{kim2021,
   author = {Kim, Jae Kyoung},
   title = {Tick, Tock, Circadian Clocks},
   booktitle = {Case Studies in Systems Biology},
   publisher = {Springer},
   address = {Cham},
   pages = {79-94},
   ISBN = {978-3-030-67742-8},
   year = {2021},
   type = {Book Section},
   editor = {Kraikivski, Pavel}
}

@article{kimtyson,
   author = {Kim, J. K. and Tyson, J. J.},
   title = {Misuse of the {M}ichaelis–{M}enten rate law for protein interaction networks and its remedy},
   journal = {{PLoS Comput. Biol.}},
   volume = {16},
   number = {10},
   year = {2020},
}

@article{feliukruffwalcher,
   author = {Feliu, Elisenda and Kruff, Niclas and Walcher, Sebastian},
   title = {{Tikhonov–Fenichel Reduction for Parameterized Critical Manifolds with Applications to Chemical Reaction Networks}},
   journal = {J. Nonlinear Sci.},
   volume = {30},
   number = {4},
   pages = {1355-1380},
   ISSN = {1432-1467},
   year = {2020},
}

@article{cabre20031,
   author = {Cabré, Xavier and Fontich, Ernest and De la Llave, Rafael},
   title = {The parameterization method for invariant manifolds I: manifolds associated to non-resonant subspaces},
   journal = {Indiana Univ. Math. J. },
   pages = {283-328},
   ISSN = {0022-2518},
   year = {2003},
}

@article{cabre20032,
   author = {Cabré, Xavier and Fontich, Ernest and De la Llave, Rafael},
   title = {The parameterization method for invariant manifolds II: regularity with respect to parameters},
   journal = {Indiana Univ. Math. J. },
   pages = {329-360},
   ISSN = {0022-2518},
   year = {2003},
}

@article{cabre2005,
   author = {Cabré, Xavier and Fontich, Ernest and De La Llave, Rafael},
   title = {The parameterization method for invariant manifolds III: overview and applications},
   journal = {J. Differ. Equ.},
   volume = {218},
   number = {2},
   pages = {444-515},
   ISSN = {0022-0396},
   year = {2005},
}

@article{multiple,
  title={Multiple timescales and the parametrisation method in geometric singular perturbation theory},
  author={Lizarraga, I. and Rink, B. and Wechselberger, M.},
  journal={Nonlinearity},
  volume={34},
  number={6},
  pages={4163},
  year={2021},
  publisher={IOP Publishing}
}

@article{mmvalid,
title = {A new {M}ichaelis-{M}enten equation valid everywhere multi-scale dynamics prevails},
journal = {Math. Biosci.},
volume = {315},
pages = {108220},
year = {2019},
issn = {0025-5564},
author = {Patsatzis, D. G. and Goussis, D. A.},
}

@article{algorithmic2023,
   author = {Patsatzis, D. G. and Goussis, D. A.},
   title = {Algorithmic criteria for the validity of quasi-steady state and partial equilibrium models: the Michaelis-Menten reaction mechanism},
   journal = {J. Math. Biol.},
   volume = {87},
   number = {2},
   pages = {27},
   ISSN = {0303-6812},
   year = {2023},
   type = {Journal Article}
}

@article{lamgoussis1989,
   author = {Lam, S. H. and Goussis, D. A.},
   title = {{Understanding complex chemical kinetics with computational singular perturbation}},
   journal = {Proc. Comb. Inst.},
   volume = {22},
   number = {1},
   pages = {931-941},
   ISSN = {0082-0784},
   year = {1989},
   type = {Journal Article}
}

@article{lamgoussis1994,
   author = {Lam, SH and Goussis, DA},
   title = {The {CSP} method for simplifying kinetics},
   journal = {Intl. J. of Chem. Kinetics},
   volume = {26},
   number = {4},
   pages = {461-486},
   ISSN = {0538-8066},
   year = {1994},
   type = {Journal Article}
}

@article{computational,
  title={Computational singular perturbation method for nonstandard slow-fast systems},
  author={Lizarraga, I. and Wechselberger, M.},
  journal={SIAM J. Appl. Dyn. Syst.},
  volume={19},
  number={2},
  pages={994--1028},
  year={2020},
  publisher={SIAM}
}

@article{buckingham,
   author = {Buckingham, Edgar},
   title = {On physically similar systems; illustrations of the use of dimensional equations},
   journal = {Phys. Rev.},
   volume = {4},
   number = {4},
   pages = {345},
   year = {1914},
   type = {Journal Article}
}

@book{foundations,
  title={{Foundations of Chemical Reaction Network Theory}},
  author={Feinberg, M.},
  year={2019},
  publisher={Springer},
series = {Applied Mathematical Sciences},
volume = {202}
}

@book{dumortierroussarie1996,
   author = {Dumortier, Freddy and Roussarie, Robert H},
   title = {Canard cycles and center manifolds},
   publisher = {Mem. Amer. Math. Soc.},
   volume = {577},
   ISBN = {082180443X},
   year = {1996},
   type = {Book}
}

@article{kimforger,
   author = {Kim, Jae Kyoung and Forger, Daniel B},
   title = {A mechanism for robust circadian timekeeping via stoichiometric balance},
   journal = {Mol. Syst. Biol.},
   volume = {8},
   number = {1},
   pages = {630},
   ISSN = {1744-4292},
   year = {2012},
   type = {Journal Article}
}

@article{gene1,
   author = {Barbuti, Roberto and Gori, Roberta and Milazzo, Paolo and Nasti, Lucia},
   title = {A survey of gene regulatory networks modelling methods: from differential equations, to {B}oolean and qualitative bioinspired models},
   journal = {J. Membr. Comput.},
   volume = {2},
   number = {3},
   pages = {207-226},
   ISSN = {2523-8906},
   year = {2020},
   type = {Journal Article}
}

@book{gene2,
   author = {Wilkinson, Darren J},
   title = {{Stochastic Modelling for Systems Biology}},
   publisher = {Chapman and Hall/CRC},
   ISBN = {1351000918},
   year = {2018},
   type = {Book}}

@book{gene3,
   author = {Bower, James M and Bolouri, Hamid},
   title = {{Computational Modeling of Genetic and Biochemical Networks}},
   publisher = {MIT press},
   ISBN = {0262524236},
   year = {2001},
   type = {Book}
}

@article{leeothmer,
  title={{A multi-time-scale analysis of chemical reaction networks: I. Deterministic systems}},
  author={Lee, Chang Hyeong and Othmer, Hans G},
  journal={	J. Math. Biol.},
  volume={60},
  number={3},
  pages={387--450},
  year={2010},
  publisher={Springer}
}

@book{kuehn2015,
  title={{Multiple Time Scale Dynamics}},
  author={Kuehn, C.},
    series = {Applied Mathematical Sciences},
  volume={191},
  year={2015},
  publisher={Springer}
}

@article{schnellmaini2000,
   author = {Schnell, S. and Maini, P. K.},
   title = {Enzyme kinetics at high enzyme concentration},
   journal = {Bull. Math. Bio.},
   volume = {62},
   number = {3},
   pages = {483-499},
   ISSN = {1522-9602},
   year = {2000},
   type = {Journal Article}
}

@article{noethenwalcher,
   author = {Noethen, Lena and Walcher, Sebastian},
   title = {Tikhonov's theorem and quasi-steady state},
   journal = {Discrete Contin. Dyn. Syst. Ser. B},
   volume = {16},
   number = {3},
   pages = {945-961},
   ISSN = {1531-3492},
   year = {2011},
   type = {Journal Article}
}

@article{schnellmaini2002,
   author = {Schnell, S. and Maini, P. K.},
   title = {Enzyme kinetics far from the standard quasi-steady-state and equilibrium approximations},
   journal = {Math. Comput. Model.},
   volume = {35},
   number = {1-2},
   pages = {137-144},
   ISSN = {0895-7177},
   year = {2002},
   type = {Journal Article}
}

@book{physio,
  title={Mathematical Physiology},
  author={Keener, J.P. and Sneyd, J.},
  year={2025},
  publisher={Springer},
series = {Interdisciplinary Applied Mathematics},
volume = {8},
edition = {3}
}

@book{param,
  title={{The
Parameterization
Method for
Invariant Manifolds}},
  author={Haro, \`A. and Canadell, M. and Figueras, J.L., Luque and A. and Mondelo, J.M.},
  year={2025},
  publisher={Springer},
series = {Applied Mathematical Sciences},
volume = {Volume 195},
}

@article{zagariskaperkaper,
   author = {Zagaris, A. and Kaper, H. G. and Kaper, T. J.},
   title = {{Analysis of the Computational Singular Perturbation Reduction Method for Chemical Kinetics}},
   journal = {J. Nonlinear Sci.},
   volume = {14},
   number = {1},
   pages = {59-91},
   ISSN = {1432-1467},
   year = {2004},
   type = {Journal Article}
}

@article{reichselkov,
   author = {Reich, JG and Sel'Kov, EE},
   title = {Mathematical analysis of metabolic networks},
   journal = {FEBS lett.},
   volume = {40},
   number = {S1},
   pages = {S112-S118},
   ISSN = {0014-5793},
   year = {1974},
   type = {Journal Article}
}

@article{laxwalcher,
   author = {Lax, Christian and Walcher, Sebastian},
   title = {Singular perturbations and scaling},
   journal = {Discrete Contin. Dyn. Syst. Ser. B},
   volume = {25},
   number = {1},
   pages = {1-29},
   ISSN = {1531-3492},
   year = {2020},
   type = {Journal Article}
}

@article{kimjosicbennett,
   author = {Kim, Jae Kyoung and Josić, Krešimir and Bennett, Matthew R.},
   title = {The relationship between stochastic and deterministic quasi-steady state approximations},
   journal = {BMC Syst. Biol.},
   volume = {9},
   number = {1},
   pages = {87},
   ISSN = {1752-0509},
   year = {2015},
   type = {Journal Article}
}

@article{kimjosicbennett2,
   author = {Kim, Jae Kyoung and Josić, Krešimir and Bennett, Matthew R.},
   title = {The validity of quasi-steady-state approximations in discrete stochastic simulations},
   journal = {Biophys. J.},
   volume = {107},
   number = {3},
   pages = {783-793},
   ISSN = {0006-3495},
   year = {2014},
   type = {Journal Article}
}

@article{coulletspiegel,
   author = {Coullet, P.H. and Spiegel, Edward A.},
   title = {Amplitude equations for systems with competing instabilities},
   journal = {SIAM J. Appl. Math.},
   volume = {43},
   number = {4},
   pages = {776-821},
   ISSN = {0036-1399},
   year = {1983},
   type = {Journal Article}
}

@article{schauerandheinrich,
   author = {Schauer, M. and Heinrich, R.},
   title = {Quasi-steady-state approximation in the mathematical modeling of biochemical reaction networks},
   journal = {Math. Biosci.},
   volume = {65},
   number = {2},
   pages = {155-170},
   ISSN = {0025-5564},
   year = {1983},
   type = {Journal Article}
}

@article{goussis,
   author = {Goussis, Dimitris A},
   title = {Quasi-steady state and partial equilibrium approximations: their relation and their validity},
   journal = {Combust. Theory Model.},
   volume = {16},
   number = {5},
   pages = {869-926},
   ISSN = {1364-7830},
   year = {2012},
   type = {Journal Article}
}

@article{robertsetal,
   author = {Roberts, A.J. and MacKenzie, Tony and Bunder, Judith E.},
   title = {A dynamical systems approach to simulating macroscale spatial dynamics in multiple dimensions},
   journal = {J. Eng. Math.},
   volume = {86},
   number = {1},
   pages = {175-207},
   ISSN = {0022-0833},
   year = {2014},
   type = {Journal Article}
}

@article{huguetetal,
   author = {Huguet, Gemma and De la Llave, Rafael},
   title = {{Computation of Limit Cycles and Their Isochrons: Fast Algorithms and Their Convergence}},
   journal = {SIAM J. Appl. Dyn. Syst.},
   volume = {12},
   number = {4},
   pages = {1763-1802},
   ISSN = {1536-0040},
   year = {2013},
   type = {Journal Article}
}

@book{roberts,
   author = {Roberts, Anthony John},
   title = {{Model Emergent Dynamics in Complex Systems}},
   publisher = {SIAM},
   ISBN = {1611973554},
   year = {2014},
series = {Mathematical Modeling and Computation},
   type = {Book}
}

@book{automatic,
   author = {Griewank, Andreas and Walther, Andrea},
   title = {{Evaluating Derivatives: Principles and Techniques of Algorithmic Differentiation}},
   publisher = {SIAM},
   ISBN = {0898716594},
   year = {2008},
series = {Other Titles in Applied Mathematics},
   type = {Book}
}

@article{eilertsenstroberg,
   author = {Eilertsen, Justin and Stroberg, Wylie},
   title = {On the reduction of stochastic chemical reaction networks},
   journal = {J. Math. Biol.},
   volume = {92},
   number = {1},
   pages = {1},
   ISSN = {1432-1416},
   year = {2025},
   type = {Journal Article}
}


\begin{tcolorbox}[
    colback=black!5,          
    colframe=black!75,        
    fonttitle=\bfseries,      
    title=Appendix: A Big Picture View of the Mathematical Methods,
    boxrule=1pt,              
    arc=2mm,                  
    ]

\subsection*{The Purpose of Mathematical Models and Dimension Reduction}

A chemical reaction network (CRN) can be thought of as a detailed `circuit diagram' for a biological process. By translating this diagram into a set of mathematical equations (here, ordinary differential equations), we create a model that allows us to simulate the process on a computer, test hypotheses, and predict how the system might behave under different conditions. However, these models are often incredibly complex, involving dozens of components and interactions. This is where {\em dimension reduction} comes in. The goal is to simplify the model by focusing only on the components and interactions that govern the system's long-term behavior. By systematically removing the fast, transient details, we create a smaller, more manageable model. This reduced model is not only faster to simulate but, more importantly, it reveals the core logic of the biological circuit.

\subsection*{An Overview of Geometric Singular Perturbation Theory}

Many biological processes involve events that happen on vastly different timescales. For example, an enzyme might bind to its substrate almost instantly (a {\em fast} process), while the final product accumulates very slowly over minutes or hours (a {\em slow} process). {\em Geometric Singular Perturbation Theory (GSPT)} is a powerful mathematical framework specifically designed to analyze systems containing this mix of fast and slow events. The core idea of GSPT is to treat the fast processes as transient paths that quickly lead the system to a stable, lower-dimensional space where the slow, meaningful dynamics unfold.

The `geometric' aspect of GSPT is key. Instead of just manipulating equations algebraically (like setting a rate to zero), this approach visualizes the problem in a multi-dimensional `state space.' The fast dynamics quickly pull the system onto a lower-dimensional object within this space, much like gravity pulls water down a landscape into a riverbed. We call this object the {\em slow manifold}. GSPT provides the tools to find this manifold and describe the `flow' along it, giving us a clear, geometric picture of the system's essential, long-term behavior. The specific methods we use (coordinate-independent GSPT and the parametrization method) are modern, highly effective ways of carrying out this geometric analysis.

\subsection*{Further Reading}

For readers interested in learning more about timescale analysis in biological systems, we recommend the following resources:

\begin{itemize}
    \item Keener, J.~P., \& Sneyd, J.~(2025). \textit{Mathematical Physiology}. Springer. 
    
    This is an excellent textbook \cite{physio} in mathematical biology. It provides an excellent introduction to timescale analysis and singular perturbations in the context of biological problems.

    \item Kuehn, C.~(2015). \textit{Multiple timescale Dynamics}. Springer. 

    A comprehensive textbook \cite{kuehn2015} that covers a wide range of methods for systems with multiple timescales, including GSPT. It is a helpful resource for understanding the broader mathematical context.

\end{itemize}
\end{tcolorbox}

\begin{tcolorbox}[
    colback=black!5,          
    colframe=black!75,        
    fonttitle=\bfseries,      
    title=Appendix: A Big Picture View of the Mathematical Methods,
    boxrule=1pt,              
    arc=2mm,                  
    ]

\begin{itemize}
    \item Wechselberger, M.~(2020). \textit{Geometric singular perturbation theory beyond the standard form}. Springer. 

    For the mathematically inclined reader, this book \cite{wechselberger2020} provides a modern treatment of the geometric theory, including the coordinate-independent methods used in this work.
    \item Haro, À., Canadell, M., Figueras, J.~L., Luque, A., \& Mondelo, J.~M.~(2016). \textit{The parameterization method for invariant manifolds}. Springer.

    {This monograph \cite{param} provides an overview of the parametrization method, from rigorous theoretical results to effective computational implementation.}
\end{itemize}

\end{tcolorbox}

\section*{Declarations}

\subsection*{Authorship and Contribution} All authors have made substantial intellectual contributions to the study conception, execution, and design of the work. All authors have read and approved the final manuscript. TEFL: formal analysis and investigation, writing, review and editing. MW: conceptualization, supervision, review and editing.

\subsection*{Conflict of Interest}
The authors declare there are no conflicts of interest.

\subsection*{Data and Code Availability}
{Code for Figure \ref{T1_diagram}b and Figure \ref{KF_discrepancy} can be found in
\url{https://doi.org/10.5281/zenodo.17717405}.}

\subsection*{Funding} TEFL acknowledges the support of an Australian Government Research Training Program (RTP) Scholarship.
MW acknowledges the support through the Australian Research Council (ARC) Discovery Project grant scheme 
(DP220101817).
\end{document}